\magnification1200
\def \sn{{\smallskip \noindent}}
\def \bn {{\bigskip \noindent}}
\def \La {{\longrightarrow}}
\def \O{{\cal O}}
\def \P{{\bf  P}}
\def \X{{\cal X}}
\def \XX{{X_{t_{\nu}}}}
\def \Xx{{X_{t_{\nu_0}}}}
\def \YY{{Y_{t_{\nu}}}}
\def \Yy{{Y_{t_{\nu_0}}}}
\def \C{{\cal C}}
\def \D{{\cal D}}
\def \Y{{\cal Y}}
\def \Z{{\cal Z}}

\def \ZZ{{Z_{t_{\nu}}}}
\def \D{{\cal D}}
\def \B{{\cal B}}
\def \Pp{{\cal P}}
\def \d{{\partial}}
\def \dd{{\overline {\partial}}}
\def \H {{\cal H}}
\def \A{{\cal A}}
\def \Pp{{\cal P}}
\def \L{{\cal L}}
\font\medium=cmbx10 scaled \magstep1
\font\large=cmbx10 scaled \magstep2

\par  \vskip.5cm
{\centerline{{\large Towards a Mori theory on compact K{\"a}hler threefolds, II}}}
\vskip.5cm
{\centerline {\medium Thomas Peternell}}
\vskip.5cm
\vskip.5cm

\noindent
{\medium  Introduction}
\sn This is the second part to my joint paper : " Towards a Mori theory on 
compact K\"ahler manifolds, 1 "  with F. Campana. With the techniques developped
in that paper we shall prove here the following result. 

\bn {\bf Main Theorem} \  {\it Let $X$ be a non-algebraic compact K\"ahler threefold satisfying
one of the following conditions.
{\item {(I)} $X$ can be approximated algebraically.
\item {(II)} $\kappa (X) = 2.$
\item {(III)} $X$ has a good minimal model.}
\sn  Assume that $K_X$ is not nef. Then 
{\item (1)} $X$ contains a rational curve C with $K_X\cdot C < 0;$
{{\item (2)} There exists a surjective holomorphic map $\varphi : X \La Y$ to a normal
complex space $Y$ with $\varphi_*(\O_X) = \O_Y$ of one of the following types.
{\item\item (a)} $\varphi $ is a $\P_1$- bundle or a conic bundle over
a non-algebraic surface (this can happen only in case (1))
{\item\item (b)} $\varphi$ is bimeromorphic contracting an irreducible divisor
$E$ to a point, and $E$ together with its normal bundle N is one of the
following $$(\P_2,\O(-1)), (\P_2,\O(-2)), (\P_1 \times \P_1, \O(-1,-1)), (Q_0,\O(-1)),$$
where $Q_0$ is the quadric cone. 
{\item\item (c)} $Y$ is smooth and $\varphi$ is the blow-up of $Y$ along a
smooth curve.}
\sn $\varphi $ is called an extremal contraction. 
\sn $Y$ is (a possibly singular) Kaehler space in all cases except possibly
(2c). Moreover in all cases but possibly (2c), $\varphi$ is the contraction
of an extremal ray in the cone ${\overline {NE}}(X).$ }  

 \bn \rm This is therefore a non-algebraic analogue to Mori theory in dimension 3.
Of course one expects that in case (2c) we can arrange things such that $Y$
is Kaehler, too, and that $\varphi$ is the contraction of an extremal ray. In principle
we would of course like to prove the theorem without any of the assumptions (I),(II) or (III).

\sn We will now explain the theorem and the method of the proof. Nefness of a
line bundle $L$ on an arbitrary compact manifold $X$ is defined via metrics :
$L$ is nef, if for every positive $\epsilon$ there is a metric $h_{\epsilon}$ on $L$
whose curvature satisfies $\Theta_{L,h_{\epsilon}} \geq - \epsilon \omega$, for a fixed
positive (1,1)-form $\omega.$ Passing to the special case $L = K_X$, it is not at 
all clear that there is any curve $C$ with $K_X \cdot C < 0,$ for arbitrary
$L$ this is even false. In order to circumvent this difficulty, we introduce the
additional alternative assumptions (I),(II) and (III).
We first focus on case (I), i.e. $X$ admits an algebraic approximation, This means that there is
a family of compact K\"ahler manifolds $(X_t)$ over the unit disc in some ${\bf C}^m$
 with $X_0 \simeq X$ and with
a sequence $t_{\nu}$ converging to 0 such that all $X_{t_{\nu}}$ are
projective. Assume now that $K_X$ is not nef. Then we will prove that 
$K_{X_{t_{\nu}}}$ is not nef for large $\nu.$ Therefore a fixed $X_{t_{\nu}}$
admits a contraction of an extremal ray. This extremal ray is given by a
non-splitting family of rational curves. The family can be deformed to the
nearby fibers and hence also to $X_0$, however the limit family may split.
Therefore one has to extract from the limit family a non-splitting family
which in turn by Part 1 of this paper will define the map $\varphi$ we are looking for. 
However it is now not clear whether $\varphi$ is induced by contractions
of extremal rays in the nearby (projective ) fibers. This causes the difficulty
in case (2c) together with some projective problem in the limit which will
be explained in full detail in sect. 4; see also below for more comments.  
\sn The proof that the bundles $K_{X_{t_{\nu}}}$ are not nef is done by arguing by contradiction.
If all $K_{X_{t_{\nu}}}$ would be nef, then a suitable multiple would be 
generated by global sections by Miyaoka-Kawamata's abundance theorem. Now we
examine the limit linear system. It turns out that in case of the presence of
a base locus, we can find a curve $C$, in the base locus, with 
$K_X \cdot C < 0$. Here we make heavily use of the fact that we are working
in dimension three (also abundance works at the moment only in this dimension)
because we have to analyse line bundles which are not nef in the analytic sense
but which are nef in the algebraic sense on non-algebraic (possibly singular)
surfaces. This is possible since the structure of surfaces of algebraic dimension $\leq 1$ is not
complicated. 
\sn If $\kappa (X) = 2$ and $K_X$ is not nef, we can directly show that there is a
rational curve $C$ with $K_X \cdot C < 0,$ using the meromorphic map (Iitaka reduction) 
attached to the linear system $\vert mK_X \vert$ for large $m.$ If $\kappa (X)  = 1,$
we can at least prove that there is a curve $C$ with $K_X \cdot C < 0.$ Once we have
a rational curve $C$ with $K_X \cdot C < 0,$ we obtain a non-splitting family of rational
curves with the same property and if $\kappa (X) 
\geq 0,$ we can apply the main result of [CP94] to obtain an extremal contraction.

\sn An important point for further developments is of course the question whether $Y$ is
again K\"ahler in case $\varphi$ is birational. For the notion of a singular K\"ahler space
we refer to sect.3. As stated already in the Main Theorem, we are able to prove this in case
that the exceptional divisor $E$ is contracted to a point, using various techiques of the
theory currents. In case ${\rm dim}\varphi (E) = 1,$ this is however not the case; at least 
there is no a priori reason why that should be true. Therefore this case remains open but
one would expect that $Y$ is K\"ahler if we have chosen the "correct" $\varphi.$ The "correct"
$\varphi$ should be given by an extremal ray in the cone of effective curves or the dual cone to
the K\"ahler cone.     

\sn Most parts of this paper have been worked out during a stay at the MSRI in
Berkeley. I would like to thank the institute for its hospitality and the
excellent working conditions.

\vfill \eject \noindent {\medium A Motivation}
\bn The final aim of a "Mori theory" of compact K\"ahler threefolds $X$ would of course be
\sn (a) to construct minimal models unless $X$ is uniruled

\sn (b) to prove abundance for minimal models: if $X'$ is a compact K\"ahler threefold with
at most terminal singularities and $K_{X'}$ nef, then $mK_{X'}$ is generated by global 
sections for suitable large $m.$ 

\bn We will discuss these topics in section 6. Besides its independent interest, the solution
of these problems would give some new insight in the structure of K\"ahler threefolds
which is not connected a priori to Mori theory. The statement we want to deduce is
the following

\bn {\it If the problems (a) and (b) have a positive solution, then simple K\"ahler threefolds
are Kummer.}

\bn The relevant definitions are :
\sn (1) a compact K\"ahler manifold is {\it simple}, if there is no covering family of positive
dimensional subvarieties (hence through a very general point of $X$ there is no positive
dimensional irreducible subvariety);
\sn (2) a compact K\"ahler manifold is Kummer, if it is bimeromorphic to a variety $T / G,$
where $T$ is a torus and $G$ a finite group acting on $T.$ 

\bn {\it Proof.} Since $X$ is simple, it cannot be uniruled. Hence by (a) $X$ has a minimal
model $X'.$ By (b), $mK_{X'}$ is generated by global sections for some $m.$ Hence the Kodaira
dimension $\kappa (X') \geq 0.$ Since $X'$ is simple, we conclude 
$ \kappa (X') = 0,$ hence
$$ mK_{X'}Ê= \O_{X'}.$$ 
Now there exists a covering
$$  \tilde X \La X',$$
unramified over the regular part ${\rm reg}ÊX,$ such that $K_{\tilde X}Ê= \O_{\tilde X},$ in particular
$\tilde X$ is Gorenstein, see e.g. [KMM87]. By the Riemann-Roch theorem for Gorenstein
threefolds (cp. e.g.[Ka86]), we obtain 
$$ \chi(\tilde X,\O_{\tilde X}) = 0. \eqno (*)$$
Now let $\hat X \La \tilde X$ be a desingularisation. Since $\hat X$ is non-algebraic we
have by Kodaira's theorem $H^2(\hat X,\O_{\hat X}) \ne 0.$ Because $\tilde X$ has only
rational singularities, we obtain
$$ H^2(\tilde X,\O_{\tilde X}) \ne 0.$$
Since $h^3(\tilde X,\O_{\tilde X}) = 1,$ we obtain from (*):
$$ H^1(\tilde X,\O_{\tilde X}) \ne 0.$$
Again working with $\hat X$ we therefore have a non-trivial Albanese map
$$ \alpha : \tilde X \La {\rm Alb}(\tilde X) = {\rm Alb}(\hat X).$$
Since $\tilde X$ is simple, $\alpha$ is surjective. Having in mind that $K_{\tilde X}Ê=
\O_{\tilde X},$ we immediately see that $\alpha$ is etale outside the singular locus
of $\tilde X.$ If $x_0$ is a singularity of $\tilde X,$ then some 1-form would
vanish at $x_0$ and hence we obtain by wedging 1-forms a 3-form with an isolated singularity
which is absurd. Hence $\tilde X$ is smooth and therefore $\tilde X$ is a torus ($\alpha$
is an isomorphism). Compare [Ka85,sect. 8]. Now $\tilde X$ being a torus, $X$ is Kummer.

\bn \bn {\medium Preliminaries }

\bn {\bf (0.1)} Let $X$ be a compact K\"ahler manifold. We say that $X$ can be approximated
algebraically, if the folowing holds. There exists a family $\X = (X_t)_{t \in \Delta}$
of compact K\"ahler manifolds over the unit disc $\Delta \subset {\bf C}^m$ such that
$X_0 \simeq X$ and there is a sequence $t_{\nu}Ê$ converging to $0$ such that all $\XX$
are projective. The projection map $\X \La \Delta$ is usually denoted $\pi.$ 
\sn By a conjecture of Kodaira (or Andreotti?) every compact K\"ahler manifolds can be
approximated algebraically, however at the moment this only known in dimension 2 (via
Enriques-Kodaira classification).

\bn {\bf (0.2)} Here we collect some standard notations. If $X$ is a compact manifold,
we let $\rho (X) $ denote the Picard number of $X$. Define $N_1(X)$ to be the linear subspace
of $H_2(X,{\bf R})$
generated by the classes of irreducible curves.  We let ${\overline {NE}}(X) \subset N_1(X)$ denote
the closed cone generated by the classes of irreducible curves.
For details on this and on Mori theory in the algebraic case in general we refer e.g. to
[KMM87].

\bn {\bf (0.3)} The algebraic dimension of an irreducible reduced compact complex space $X$ is by definition
the transcendence degree of its field ${\cal M}(X)$ of meromorphic functions over ${\bf C}.$ 
\sn An algebraic reduction of $X$ is a meromorphic map $f: X \rightharpoonup Y$ inducing an isomorphism
$f^*: {\cal M}(Y) \La {\cal M}(X).$  

\bn {\bf (0.4)} For the convenience of the reader we state here the main result of [CP94] which
we shall use several times. We shall use the notion of a non-splitting family $(C_t)_{t \in T}$
of rational curves. This means that $T$ is compact and irreducible and every $C_t$ is an
irreducible and reduced rational curve. Now the main result of [CP94] reads as follows

\bn {\bf Theorem}  {\it Let $X$ be a compact K\"ahler threefold  and $(C_t)_{t \in T}$ a
nonsplitting family of rational curves. 
{\item {(1)} If $K_X \cdot C_t = -4 $ and ${\rm dim}T = 4,$ then $X \simeq \P_3.$
\item {(2)} If $K_X \cdot C_t = -3$ and ${\rm dim}T = 3,$ then either $X \simeq Q_3,$ the
threedimensional quadric, or $X$ is a $\P_2-$bundle over a smooth curve.
\item {(3)} Assume $K_X \cdot C_t = -2 $ and ${\rm dim}T = 2.$ 
\itemitem {(3.1)} If $X$ is non-algebraic and the $(C_t)$ fill up a surface $S \subset X,$ then
$S \simeq \P_2$ with normal bundle $N_{S \vert X}Ê= \O(-1)$ (the same holds for $X$ projective
 if $S$ is normal)
\itemitem {(3.2)} If $X$ is covered by the $C_t$, then one of the following holds.
\itemitem {(3.2.1)} $X$ is Fano with $b_2(X) = 1$ and index 2.
\itemitem {(3.2.2)} $X$ is a quadric bundle over a smooth curve with $C_t$ contained in fibers
\itemitem {(3.2.3)} $X$ is a $\P_1-$bundle over a surface, the $C_t$ being the fibers.
\itemitem {(3.2.4)} $X$ is the blow-up of a $\P_2-$bundle over a curve along a section. Here
the $C_t$ are the strict transforms of the lines in the $\P_2$'s meeting the section
\item {(4)} Let $K_X \cdot C_t = -1$ and ${\rm dim}T = 1.$ Then the $C_t$ fill up a surface 
$S.$ Assume $X$ non-algebraic.
\itemitem {(4.1)} If $S$ is normal, then one of the following holds.
\itemitem {(4.1.1)} $S = \P_2$ with $N_S = \O(-2).$
\itemitem {(4.1.2)} $S = \P_1 \times \P_1$ with $N_S = \O(-1,-1).$
\itemitem {(4.1.3)} $S = Q_0$ (a quadric cone) with $N_S = \O(-1).$ 
\itemitem {(4.1.4)} $S$ is a ruled surface over a smooth curve and $X$ is the blow-up of
a smooth threefold along $C.$
\itemitem {(4.2)} Let $S$ be non-normal. Then $\kappa (X) = - \infty.$ If moreover $X$ can
be approximated by algebraic threefolds, then we have $a(X) = 1$ and $X$ has the structure of
a "generic conic bundle" over a surface $Y$ with $a(Y) = 1.$ The general fiber of the algebraic
reduction $f$ which can be taken holomorphic in this case is an almost homogeneous 
$\P_1-$bundle over an elliptic curve. The surface $S$ consists of reducible conics and is
contracted by $f$ to a point.}}

\bn The essential content of the theorem can be rephrased as follows. Assume that $C$ is a
rational curve with $K_X \cdot C = k, -1 \geq k \geq -4.$ If no deformation of $C$ splits,
then the conclusions of the theorem hold. The point is that automatically the dimension
of the deformation is at least $k.$  
 
\bn We will have a more careful look at the case (4.2) in sect.5.    

\bn The paper [CP94] will be refered to as Part I.

\bn \bn

{\noindent \medium 3. Limits of extremal rays}
\bn {\bf (3.1)} For most of this section  we fix a family $\pi : \X \La \Delta $ of compact
K\"ahler threefolds $X_t = \pi^{-1}(t)$ over the unit ball $\Delta \subset {\bf C}^m.$
Let $t_{\nu}$ be a sequence converging to 0. Assume that all $\XX$ are projective, whereas
$X_0$ is {\bf not}. In other words, $\X$ is an algebraic approximation of $X = X_0.$  
Furthermore we assume that all $K_{\XX}$ are not nef. Fix some $\nu_0$ and consider an extremal
ray $R_{\nu_0}$ on $\Xx$ generated by the extremal rational curve $C_{\nu_0}.$ Let
$\varphi_{\nu_0} : \Xx \La \Yy$ be the contraction associated with $R_{\nu_0}.$ We
are going to construct from $R_{\nu_0}$ a sequence of extremal rays $R_{\nu}$ living on $X_{\nu}$. 
In order to do this we first examine the structure of $R_{\nu_0}.$ Let $m$
denote the length of $R_{\nu_0},$ i.e.
$$ m  = {\rm min} \{ -K_{\Xx}Ê\cdot C \vert [C]Ê\in R_{\nu_0} \}.$$

\sn Then we have :

\bn {\bf 3.2 Lemma } {\it Either $m = 1$ or $m = 2.$ }

\bn {\bf Proof.} Assume $m \geq 3.$ Because of the a priori bound $m \leq 4$ we have
only to exclude the cases $m = 4$ and $m = 3.$ In the first case we have $\Xx \simeq
\P_3,$ in the second $\Xx \simeq {\bf Q}_3,$ the threedimensional quadric or $\Xx$ is
a $\P_2-$ bundle over a curve [Wi89]. Hence always $H^2(\Xx,\O_{\Xx}) = 0.$ Therefore
$H^2(X_0, \O_{X_0}) = 0$ and $X_0$ is projective by Kodaira's well-known theorem,
contradiction.

\bn The same argument shows 

\bn {\bf 3.3 Lemma} {\it The contraction $\varphi_{\nu_0}$ cannot be a del Pezzo
fibration (i.e. ${\rm dim }\Yy \ne 1$) nor a $\P_1-$ or a conic bundle over a surface of Kodaira 
dimension  $- \infty.$}

\bn In fact, otherwise we would have $H^2(X_0, \O_{X_0}) = 0$ by virtue of
$$R^{i}\varphi_{\nu_0*}(\O_{\Xx}) = 0, i \geq 1, \  {\rm and} \   H^2(\Yy,\O_{\Yy}) = 0.$$ 

\bn {\bf (3.4)} From [Mo82] we obtain the following structure of $\varphi_{\nu_0}:$
{\item {(a)} a $\P_1-$bundle over a surface with nonnegative Kodaira dimension
\item {(b)} a conic bundle over a surface with nonnegative Kodaira dimension
\item {(c)} a birational contraction contracting an irreducible divisor $E$ such that
\itemitem {(c.1)} $E = \P_2$ with normal bundle $N = \O(-1),$
\itemitem {(c.2)} $E = \P_2, N = \O(-2),$
\itemitem {(c.3)} $E = \P_1 \times \P_1, N = \O(-1,-1),$
\itemitem {(c.4)} $E = Q_0,$ the quadric cone, with $N = \O(-1).$ 
\item {(d)} a birational contraction contracting an irreducible divisor $E$ such that
 ${\rm dim}\varphi_{\nu_0} = 1 $ and $\varphi_{\nu_0}$ is the
blow-up of the smooth curve $\varphi_{\nu_0}(E)$ in the manifold $\Yy.$} 

\sn We will call the different contractions of type (a), (b), (c.i) and (d), respectively.

\bn Now consider the extremal rational curve $C_{\nu_0}Ê\subset \Xx.$ We assume that
$-K_{\Xx}Ê\cdot C_{\nu_0} $ is minimal, $C_{\nu_0}$ is then smooth. Then the normal bundle 
$N$ of $C_{\nu_0}$ in $X_{t_{\nu_0}}$ is of the form 
$$ N = \O(a) \oplus \O(b) $$
with $(a,b) \in \{(0,0),(0,-1),(1,-1),(1,-2)\}.$ 
\sn If $a,b \geq -1,$ the deformations of $C_{\nu_0}$ in $\Xx$ as well as in $\X$ are
unobstructed, since 
$$ N_{C_{\nu_0} \vert \X} \simeq \O(a) \oplus \O(b) \oplus \O.$$ 
Let $\C = (C_s)_{s \in S}$ be the family of deformations of $C_{\nu_0}$ in $\X,$ in particular we
obtain a limit family $(C^0_s)$  in $X_0.$ However it is not clear whether e.g. a
$\P_1-$bundle structure on $\Xx$ converges to a $\P_1-$bundle structure on $X_0.$
More precisely, the subfamily $(C_s)_{s \in S_{\nu_0}}$ of deformations of $C_{\nu_0}$
inside $\Xx$ forms a non-splitting family of rational curves in the sense of Part 1;
however the limit family may split. The simplest example (in the case of projective
families) is the specialisation of $\P_1 \times \P_1$ into the Hirzebruch surface
$\P(\O \oplus \O(-2)) $  where one of the two rulings converges to a splitting family
of rational curves. 
\sn We also want to have a limit family in case $(a,b) = (1,-2).$ This case occurs
if either $\varphi_{\nu_0}$ is a conic bundle, here $C_{\nu_0}$ must be the reduction
of a non-reduced conic, or $E = \P_2$ with $N = \O(-2).$ In the first case we just take
$C_{\nu_0} $ not to be the reduction of a non-reduced conic, but an irreducible 
component of a reducible reduced conic; then $N = \O \oplus \O(-1)$ and we can conclude.
If $E = \P_2$ with $N = \O(-2)$ then by deformation theory the deformations of $E$ in
$\X$ are unobstructed, so we can deform $E$ in $\X$ in a 1-dimensional family (note
$N_{E \vert \X} = \O(-2) \oplus \O$). So our family $\C$ exists also here. 

\bn {\bf 3.5 Theorem} {\it There exists a surjective holomorphic map $$\psi_0 : X_0 \La Z_0$$
to a normal complex space $Z_0$ in class $\C$ of type (a),(b), (c.1) or (c.2) , as 
described in (3.4), such that $\rho(X) = \rho(Y) + 1.$ Moreover there exists a morphism
$\psi_t : X_t \La Z_t$ for general $t$ of the same type fitting into a family 
$\psi : \X \La \Z$ outside a proper analytic set $A \subset \Delta$ not containing $0.$ 
The $\ZZ$ are projective except possibly in case (c.2).}

\bn Recall that $\rho(X) $ is the Picard number of $X$ and that a normal compact 
complex space is in class $\C$ if it is bimeromorphically equivalent to a K\"ahler
manifold. 

\bn {\bf Proof.} We consider the family $\C$ constructed in (3.4). In particular we have
the family $(C^0_s)$ in $X_0.$ For simplicity of notations we skip the upper index $0.$ 
\sn (1) Assume first that the family $(C_s)$ does not split. Then we apply the Main Theorem
of Part 1 (0.4) and obtain the following. \sn If $-K_{X_0}Ê\cdot C_s = -K_{\Xx} \cdot C_{\nu_0} = 2,$
then, $X_0$ being non-algebraic, the $C_s$ define either a $\P_1$-bundle structure or
a birational contraction of type (c.1) with $E = \P_2$ and $N = \O(-1).$ This defines
our morphism $\psi_0 : X_0 \La Z_0.$ It is clear that the corresponding families in the nearby 
fibers defines maps $\psi_t$ of the same type and that all maps fit together. Moreover
$\psi_{t_{\nu_0}} = \varphi_{t_{\nu_0}}.$ 
\sn Now let $-K_{X_0} \cdot C_s = 1.$ Let $E_0 = \bigcup C_s \subset X_0.$ If $E_0$ is
normal, then we are in case (c.i) ($2 \leq i \leq 4$) or (d) and everything holds in the same way as above. So
let $E_0$ be non-normal. Then $X_0$ has the structure of a so-called generic conic bundle
in the sense of Part 1. If the $C_s$ define an extremal ray in ${\overline {NE}}(X_0)$, or,
equivalently, if the associated full family of "conics" splits only in the standard way
in the sense of (5.1), then we conclude by (5.2) resp. (5.3). Otherwise we consider the
general smooth conic and the associated family $(\tilde C_t).$ Then we have a non-standard
splitting
$$ \tilde C_{t_0} = \sum_{i=1}^p C'_i $$
with $p \geq 3$ or with $p = 2$ and, say, $-K_{X_0}Ê\cdot C'_1 \geq 2.$ In any case we have
a new $C'_i$ with $-K_{X_0} \cdot C'_i \geq 1$ which deforms in $X_0$. Then we proceed 
as in (2) below. 
\sn (2) Assume now that $(C_s)$ splits :
$$ C_{s_0} = \sum C'_i.$$
Say that $-K_{X_0} \cdot C'_1 > 0.$ Then $C'_1$ deforms in an at least 1-dimensional family,
say $(C'_t)$, in $X_0.$ If this family does not split, then we argue as in (1); note that
the $C'_t$ deform also to the nearby fibers since in the normal bundle we add a trivial
factor (see (3.4) for details). If the family splits again, we take a splitting part $C''_1$
with $-K_{X_0} \cdot C''_1 > 0$ etc. The only thing left is to prove that this procedure
terminates. 
\sn (a) First assume that no generic conic bundle comes up in our procedure. Then
we consider the cone ${\overline {NE}}(X_0) \subset N_1(X_0);$ note that every effective
curve represents an integer point in ${\overline {NE}}(X_0),$ denote this set of integer
points by ${\overline {NE}}(X_0)_{\bf Z}.$ Now introduce a norm $\Vert \cdot \Vert $
in $N_1(X_0) $ such that $\Vert x \Vert ^2 \in {\bf N}$ for $x \in {\overline {NE}}(X_0)_{\bf Z}.$ 
Let $x_1$ be the class of a general element of our first family $(C_s)$, $x_2$ the one of
the second etc. \sn Then $x_1 = x_2 + x_2' $ with $0 \ne x_2 \in {\overline {NE}}(X_0)_{\bf Z}.$ 
Therefore $$\Vert x_2 \Vert < \Vert x_1 \Vert $$ (note that ${\overline {NE}}(X_0) \cap
-{\overline {NE}}(X_0) = \{0 \}$ since $X_0$ is K\"ahler). After finitely many steps we
reach $\Vert x_n \Vert = 0,$ so $x_n  = 0,$ which means that the family in step $n-1$
does not split. 
\sn (b) If generic conic bundles come into to the game, it seems at first glance that 
they cause difficulties because we then pass from $x$ to $2x$. However this difficulty disappears 
 if $X_0$ carries only finitely many generic conic bundle structures. Indeed, $X_0$
can carry at most one generic conic bundle structure $X_0 \La Z_0$, because a second
one would produce another covering family of rational curves on $Z_0$ so that $Z_0$ would be
algebraic, hence $X_0$ would be algebraic. 
\sn The claim on the Picard number is clear, since by construction a curve $C$ is
contracted by $\psi_0$ iff $C \equiv a C_s.$
\sn Finally, $\psi_0$ extends to $\psi_t$ for $t$ near $0,$ because our non-splitting family
of rational giving rise to $\psi_0$ extends to non-splitting families on the $X_t$, therefore
the claim follows rather easily from (0.4).  

\bn {\bf 3.6 Addendum} \ Ê{\it In cases (a), (b),(c.1) and (d) the space $Z_0$ is smooth. 
In particular $Z_0$ is
a K\"ahler surface in (a) and (b), and the curves contracted by $\psi_0$ form an
extremal ray in ${\overline {NE}}(X_0).$ In the case of (c.1) it is classical
that $Z_0$ is K\"ahler, too, hence we have the same conclusion.} 

\bn It is clearly important to know whether $Z_0$ is K\"ahler also in the remaining cases.

\bn { \bf 3.7 Theorem} \  {\it Let $X$ be a (smooth) compact K\"ahler threefold, 
$\varphi : X \La Y $ be a birational morphism contracting the irreducible divisor 
$E$ to a point $p.$ Let $N_E$ be its normal bundle. Assume that $(E,N_E)$ is one of
the following : \sn $(\P_2, \O(-2)), (\P_1 \times \P_1, \O(-1,-1)), (Q_0, \O(-1)).$ 
In case $E = \P_1 \times \P_1,$ assume also that $s \times \P_1 \equiv \P_1 \times t
$ in $X_0$ for all $s,t.$  \sn Then $Y$ is a normal K\"ahler space.}  

\bn {\bf 3.8 Explanations} \  (1) A reduced complex space $X$ is K\"ahler (cp. Grauert [Gr62])
if there is a K\"ahler form $\omega$ on the smooth part reg$X$ such that the following
holds. For every singular point $x \in X$ there exists an open neighborhood $U \subset X,$
an open set $V \subset {\bf C}^N,$ a biholomorphic map $f : U \La A \subset V$ onto a closed 
subspace A and a K\"ahler form $\omega'$ on $V$ with $$f^*(\omega') \vert {\rm reg}X \cap U
= \omega .$$ 
\sn (2) Let us explain (3.7) in case where $X$ is projective. Then $Y$ is projective, too. 
In fact, let $L$ be an ample line bundle on $X.$ Then we can write
$L \vert E = \O_E(-kE)$ with some $k > 0$ (in case $E = \P_2$ we eventually substitute
$L$ by $L^{\otimes 2}).$ Note that $k$ exists in case $E = \P_1 \times \P_1$ because of
our assumption $s \times \P_1 \equiv \P_1 \times t.$ Now there exists a line bundle $L'$
on $Y$ with 
$$ L \otimes \O_X(kE) \simeq \varphi^*(L'),$$
and it is immediately clear that $L'$ is ample .         
\sn (3) Assume $X$ projective and $E = \P_1 \times \P_1  $ with $s \times \P_1 \not\equiv 
\P_1 \times t.$ We want to describe the situation; this is of course well known. 
 Consider the blow-down $\varphi_i : X \La Y_i$ contracting $E$ along 
the two projections to curves $C_i \subset Y_i.$ Then $C_i \simeq \P_1$ with normal bundle
$N_{C_i} = \O(-1) \oplus \O(-1). $ So we obtain birational maps $\psi_i : Y_i \La Y$
contracting exactly the $C_i.$ Then $Y$ is projective if and only if $\varphi $ is the
contraction of an extremal face which just means that the face generated by the fibers of
$\varphi_i, i=1,2$ is extremal. 
\sn On the other hand there are examples where $Y$ is not projective. 

\bn {\bf Proof of 3.7} (I) In a first step we construct a semi-positive closed (1,1)-form
$\omega'$ on $X$ which is positive on $X \setminus E$ such that $\omega' \vert E = 0.$
 \sn (I.1) Fix a K\"ahler form $\omega $ on $X.$ 
As in the projective case (3.8(2)) we can choose $\lambda > 0$ such that 
$$ [\omega] + \lambda c_1(\O_X(E)) \vert E = 0 \eqno (*) $$
in $H^2(E,{\bf R}).$ \sn  Fix a representative $\eta$ of $c_1(\O_X(E)).$
We note that if $\eta'$ is another representative of $c_1(\O_X(E)),$ then we have 
$\eta - \eta' = \dd \rho ,$ but we have even more, namely
$$ \eta - \eta' = \d \dd h  \eqno (+)$$
with a $C^{\infty}-$function $h$ on a neighborhood $U$ of $E.$ This is completely standard:
\sn since $\dd \d \rho = 0,$ $ \partial \rho $ is a holomorphic 2-form, therefore by the holomorphic $d-$
Poincar\'e lemma (note that $H^1(U,{\bf C}) = H^1(E,{\bf C}) = 0$ choosing $U$ such that  $E$ is
a retract of $U$) we can write $\d \rho = d \varphi$ with a holomorphic 1-form $\varphi$. 
Now put $ \tilde \rho = \rho - \varphi;$ then $\dd \tilde \rho = \eta - \eta'$ and 
$\d \tilde  \rho = 0.$ Hence $ \tilde \rho = \d h$ by Dolbeault (note $H^1(U,\O_U) = 0$ if $U$ is
strongly pseudo-convex) and (+) follows.

\bn (I.2) Let $\D$ denote the space of bidimension (1,1)-currents on $X.$ Let $W \subset \D$
be the subspace of those $\partial \overline {\partial}$  closed currents with ${\rm supp}T \subset E.$ 
Obviously  $W$ is closed. Following [HL83]
we let $\Pp \subset \D$ be the  closure of the cone of the positive (1,1)-currents $T$. Let $\B \subset \D$ be  
the subspace given by the $(1,1)$-parts of the currents of the form $dS.$ Letting 
$\kappa : \D \La \D / W =
\tilde\D$ be the projection, we set $\tilde \Pp = \kappa(\Pp), \tilde \B = \kappa (\B).$
Then we claim :
$$ \tilde \PpÊ\cap \tilde \B = \{0\}.  \eqno (**) $$ 
In fact, let $0 \ne T \in \D$ with $\kappa (T) \in \tilde \Pp \cap \tilde  \B.$ We may assume that $T \in 
\Pp.$ Then there are currents $T'$ and $S$ with supp$T' \subset E, dT' = 0$ such that 
$$ T = T' + (dS)_{(1,1)}.  $$
It follows that $\partial {\overline \partial}T = 0.$

\sn Decompose $T = \chi_ET + \chi_{X \setminus E}T.$  A difficulty arising is that 
$T$ is in general only $\d \dd-$closed and not
$d-$closed. Let us first prove (**) in the simpler case $dT = 0$. Then by [Sk82] we have
$d(\chi_ET) = 0.$ Now observe that $T'(\omega + \lambda \eta) = 0$ by (+) because we can choose
$\eta $ in such a way that $\omega + \lambda \eta = 0$ in a neighborhood of $E$ (see I.3). Hence we get
$$ T(\omega + \lambda \eta ) = 0.$$
This is still independent on the closedness assumption.
By the same reason as for $T'$, we also have $\chi_ET(\omega + \lambda \eta) = 0.$ Here we have
used $d\chi_ET = 0,$ so that we can take $\eta $ as we want. 
Hence $\chi_{X \setminus E}T(\omega + \lambda \eta) = 0.$ On the other hand, $\chi_{X \setminus E}T
(\omega) \geq 0$ and by (3.7.a), also $\chi_{X \setminus E}T(\eta) \geq 0,$ because this does
not depend on the choice of $\eta$ ! Therefore $\chi_{X \setminus E}T = 0$ which was to be proved.
\sn The difficulty in the general case is to prove $$\d \dd \chi_ET = 0. \eqno (A)$$ 
Once we know (A) we can conclude as before; since (A) 
implies that $\chi_ET(\eta)$ (and hence $\chi_{X \setminus E}T(\eta))$
does not depend on the choice of $\eta.$ Clearly (A) will follow from 
$$ \d \dd \chi_ET \geq 0,Ê\eqno (B)$$
where $\geq$ is to be understood in the sense of currents (just apply (B) to a constant
non-zero function).  
\sn In order to prove (B) we make use  of the following theorem on currents.
\sn (C) {\it Let $\Omega \subset {\bf C}^n$ be an open set, $E \subset \Omega$ be a complex
submanifold and $T$ a positive current with $\d \dd T \geq 0$ (called plurisubharmonic in the
literature). Then $\d \dd \chi_ET \geq 0.$}
\sn This follows from Bassanelli's paper [Ba94,1.24,3.5], as was pointed to me by L.Alessandrini
[Al96]. Now (C) implies our claim (B) in case $E$ is smooth. In case $E$ is singular, i.e. a
quadric cone with vertex $x_0,$ we need an additional argument. By (B) we certainly know that
$$ \d \dd \chi_ET \vert X \setminus \{x_0 \}Ê$$
is positive, so does its trivial extension $(\d \dd \chi_ET)^0.$
Since $$\chi_ET = (\chi_ET \vert X \setminus \{x_0\})^0 $$ (note $\chi_{\{x_0\}}T = 0,$ [Ba94,1.13]),
we obtain from [AB93,5.11]
$$ \d \dd \chi_ET = (\d \dd \chi_ET)^0 + \lambda T_{x_0},$$
where $\lambda \leq 0$ and $T_{x_0}$ is the Dirac distribution of $x_0.$ 
We want to prove $\lambda = 0.$ Take a smooth closed (2,2)-form $u \leq 0$ such that the
current $T + u$ is negative: $$T + u \leq 0.$$ E.g. consider the 
K\"ahler form $\omega $ and let $u = -k \omega \wedge \omega $ for some large $k.$ Having in
mind $\chi_E(T + u) = \chi_ET$ and
$$ \d \dd \chi_E(T + u) \geq 0$$
on $X \setminus \{x_0\},$ we can apply [AB93] and obtain 
$$ \d \dd \chi_E(T + u) = (\d \d \chi_E(T + u) \vert X \setminus \{x_0 \})^0 + \mu T_{x_0}$$
with $\mu \geq 0.$
Therefore we conclude $\lambda = \mu = 0$ and we get $\d \dd \chi_ET = 0$ also in the case
of the quadric cone.

\sn (I.3) To finish the proof of (**), we finally show that we can choose a representative of $\omega
+ \lambda \eta,$ which is $0$ on $U.$ For this one needs to solve the equation
$$ \omega + \lambda \eta = - \dd \rho $$
on $U$ which leads to ask for the injectivity of
$$ H^1(U,\Omega^1_U) \La H^1(E,\Omega^1_E).$$
This follows immediately from
$$ H^1(E,N^*_E) = H^1(E,N^{*\mu}_E \otimes \Omega^1_U) = 0, \mu \geq 1.$$
Note that these arguments also work in the case $E$ is the quadric cone, we leave the 
details to the reader.

\bn (I.4)) Having verified (**) we find by the Hahn-Banach theorem a linear functional $\tilde \Phi : \tilde \D \La {\bf R}$
which is strictly positive on $\tilde \Pp \setminus 0$ and $0$ on $\tilde \B.$ Lift 
$\tilde \Phi $ to a functional $\Phi : \D \La {\bf R}.$ Then $\Phi $ is given by a $\C^{\infty} 
(1,1)-$form $\omega'$ with the following properties.
{\item {(a)} $d\omega' = 0$}
{\item {(b)}Ê$\omega' \vert X \setminus E$ is positive, so $\omega'$ itself is semipositive}
{\item {(c)} $[\omega' \vert E] = 0$ in cohomology.}

\sn (b) and (c) together give $\omega' \vert E = 0.$ 
\bn (II) We next claim that $\omega'$ is induced by a $(1,1)-$form $\omega_0$ on $Y$ and which is
almost a K\"ahler metric on $Y.$ Of course, $\omega_0$ exists on $Y \setminus p,$ where $p = \varphi(E).$
We show that there is an open neighborhood $V$ of $p$ and a plurisubharmonic function $g$
on $V$ such that $\omega_0 = \d \dd g$ on $V \setminus p.$ Then $\omega_0$ exists as form on
$Y.$ However it might not quite be a K\"ahler form on $Y$ since $g$ is possible not strictly
plurisubharmonic at $p.$ But then take a closed $(1,1)-$form $\lambda$ on $Y$ which is positive
at $p$ and let $\tilde \omega = \omega_0 + \epsilon \lambda. $ For sufficiently small $\epsilon,
\tilde \omega$ will be a K\"ahler form on $Y.$ So it remains to prove the existence of $g.$ 
For this we need to construct a $\C^{\infty}-$function $f$ on a suitable neighborhood $U$ of
$E$ such that $f \vert U \setminus E$ is strictly plurisubharmonic and $\omega' = \d \dd f $
on $U \setminus E.$ Then automatically $f = \varphi^*(g).$ 
\sn Let $\H$ denote the sheaf of pluriharmonic functions on $X$. Taking real parts of
holomorphic functions, we have an exact sequence 
$$ 0 \La {\bf R} \La \O_X \La \H \La 0 \eqno (S). $$
Let $U$ be a strongly pseudo-convex neighborhood of $E$, such that  $E$ is a deformation retract
of $U.$ By [HL83] we have 
$$ {H^1(U,\H) \simeq} {{\{Ê\psi \in \A^{1,1}_{\bf R}(U) \vert d\psi = 0\}}Ê\over {\d \dd \A^{0,0}_{\bf R}
(U)}} ,$$
where $\A^{p,q}_{\bf R}$ is the sheaf  of real $(p,q)-$forms. Therefore $\omega' \vert U$
defines a class $[\omega']Ê= 0 \in H^1(U,\H_U)$ and we must show $[\omega'] = 0.$ 
Since $\omega' \vert E = 0, $ we have $[\omega' \vert E] = 0$ in $H^1(E, \H_E)$, and so it
is sufficient to verify that  the restriction map 
$$ H^1(U, \H_U) \La H^1(E,\H_E) \eqno(R) $$
is an isomorphism. 
By the exact sequence (S) on $U$ resp. $E$ and the obvious facts 
$$H^{i}(U,\O_U) \simeq H^{i}(E,\O_E) = 0, i=1,2 $$
we have $$H^1(U,\H_U) \simeq H^2(U,{\bf R}) ,$$
  $$H^1(E,\H_E) \simeq H^2(E,{\bf R}).$$
Hence (R) follows from the fact that 
$$ {\rm rest} : H^2(U,{\bf R}) \La H^2(E,{\bf R})$$
is an isomorphism, $E$ being a deformation retract of $U.$ So $[\omega' \vert U] = 0$
in $H^1(U,\H_U),$ and $\omega' \vert  U = \d \dd f$ with  $f \in \A^{0,0}(U).$ 
Since $\omega'$ is positive on $X \setminus E$,  the function $f$ is strictly plurisubharmonic
on $U \setminus E.$

\bn {\bf 3.7.a Sub-Lemma} {\it Let $X$ be a compact K\"ahler manifold, $D$ an irreducible 
divisor on $X.$ Let $T$ be a positive closed current on $X$ of bidegree $(n-1,n-1),$
where $n = {\rm dim}X.$ Assume $\chi_DT = 0.$ Then the intersection number $T \cdot D \geq 0.$}

\bn {\bf Proof.} I am indepted to Jean-Pierre Demailly for communicating to me the following
short proof using his approximation theorem for positive closed curents [De92]. We can write $D$ (the current
of integration over $D$)
as a weak limit of smooth closed forms $\Theta_{\varepsilon}$ in the same cohomology class 
as $[D]$ such that
$$ \Theta_{\varepsilon} \geq - \lambda_{\varepsilon}Êu - O(\varepsilon) \omega,$$
where $\omega $ is a positive (1,1)-form, $u$ a suitable semi-positive (1,1)-form on $X$
depending on the global structure of $X$ and $(\lambda_{\epsilon})$ a decreasing family 
of non-negative smooth functions $(0 < \varepsilon < 1)$ converging pointwise to $0$ on
$X \setminus D$ and to the multiplicity $m(D,x)$ on $D.$  
It follows
$$ D \cdot T = \Theta_{\varepsilon} \cdot T \geq  -\int_X \lambda_{\varepsilon}Êu \wedge T - 
O(\varepsilon).$$
Now the monotone convergence theorem gives convergence to $0$, because $\chi_DT = 0.$

\bn {\bf 3.7.b Remark} \  $T\cdot D$ can be computed either by $[T]Ê\cdot [D],$ where $[T]$ and $[D]$ are the
classes in $H^2(X,{\bf R})$ resp. $H^{2n-2}(X,{\bf R})$ or by representing $D$ by a positive
closed $(1,1)-$form (the curvature form of a metric on the line bundle associated to $D$),
say $\eta,$ and computing $T(\eta).$ 
\sn If merely $\d \dd T = 0,$ then we still can define a class $[T]$, as explained detailed in
the proof of (3.11), and everything in (3.7.a) remains true (the regularisation is applied to
the current $D$ and not to $T!)$

\bn {\bf 3.7.c Remark}  Theorem 3.7 should be true in a more general context : Let $X$ be
a compact K\"ahler manifold, $f: X \La Y$ a birational morphism to a normal complex space.
Assume that the exceptional set of $f$ is an irreducible divisor $E$ and that $\rho (X) =
\rho (Y)  + 1.$ Then $Y$ is K\"ahler. 
\sn The proof should be the same as before except that one needs a singular version of
Bassanelli's theorem [Ba94,3.5] which must be true.

\bn {\bf 3.9 Corollary} {\it Let $\varphi : X \La Y$ be a contraction of type (c).
If $\rho(X) = \rho(Y) + 1,$ then $Y$ is K\"ahler and $R = {\bf R}_+[l]Ê\subset {\overline
{NE}}(X)$ is an extremal ray, where $l$ is any curve contracted by $\varphi.$ \sn
If conversely $R$ is an extremal ray, then obviously $\rho(X) = \rho(Y) +1$ and $Y$ is
K\"ahler. }

\bn {\bf Proof.} (a) If $\rho(X) = \rho(Y) +1, $ then the case $E = \P_1 \times \P_1, {\rm dim}\varphi (E) = 0, 
s \times \P_1 \not \equiv \P_1 \times t $ is excluded, hence 3.7 applies and it is 
immediately checked that $R$ is extremal. 
\sn (b) If $R$ is an extremal ray then the same arguments apply.

\bn {\bf 3.10 Remarks} (1) The case $E = \P_1 \times \P_1, {\rm dim} \varphi(E) = 0, s \times \P_1 \not \equiv
\P_1 \times t $ is not relevant for us, since then the associated contraction is
not defined by an extremal ray and not given by one non-splitting family of
rational curves, cp. (3.8).
\sn (2) Clearly the proof of (3.7) also works at least in the case where $X$ is a compact
K\"ahler manifold of dimension $n$ and where $\varphi: X \La Y$ is an extremal contraction
contracting a smooth divisor to a point.

\bn It remains to treat birational contractions of type (d), i.e. blow-ups of smooth curves.
 Already in 
 case $X_0$ is projective, $Y_0$ will not be projective "in general". This happens
exactly if the defining ray $R$ is extremal. In order to generalise this to the
K\"ahler case we consider the closed dual cone ${\overline {NA}}(X) $ to the K\"ahler cone.
So ${\overline {NA}}(X) \subset H^4(X,{\bf R})$ resp. $H^{2,2}(X).$ In general, ${\overline
{NE}}(X) $ is a proper subcone of ${\overline {NA}}(X)$ and we can view ${\overline {NA}}(X)$
as the cone generated by the classes of positive closed currents of bidegree $(2,2).$ 

\bn {\bf 3.11 Proposition} {\it Let $X$ be a compact K\"ahler threefold, $\varphi:
X \La Y$ a birational map of type (d), i.e. $X$ is the blow-up of a smooth curve
$C$ in the manifold $Y.$ Let $E = \varphi^{-1}(C)$ and $l \subset E$ a ruling line. 
Then $Y$ is K\"ahler if and only if $R = {\bf R}_+[l]$ is extremal in ${\overline {NA}}(X).$}

\bn {\bf Proof.} (a) Assume that $Y$ is K\"ahler. In order to show that  $l$ is  extremal,
we take
positive closed currents $T_1, T_2 $ with $[l] = [T_1] + [T_2] $. 
Then $0 = \varphi_*(T_1 ) + \varphi_*(T_2) ,$ and since $Y$ is K\"ahler, we have
for a K\"ahler form $\omega$ on $Y$ that $\varphi_*(T_i)(\omega) = 0,$ hence 
$\varphi_*(T_i) = 0$. Therefore the $T_i$ are supported on fibers of $\varphi$, i.e. the
$T_i$ are linear combinations of currents "integration over a fiber" and therefore
$[T_i]Ê\in R.$ 
\sn (b) Conversely, assume that $R$ is extremal in ${\overline {NA}}(X).$ By [HL83]
it is sufficient to show the following: 
\sn (*) if $T$ is a positive current of bidegree (2,2) with  $T = {\d} {\overline S} + {\dd S} ,$ 
then $T = 0.$ 
\sn For the proof of (*) write 
$$ T = \chi_CT + \chi_{Y \setminus C}T .$$
If $dT = 0,$ then by [Si74] $\chi_CT = \lambda T_C$ where $T_C$ is the current "integration over 
$C."$ In the general case this is due to Bassanelli [Ba94]. Putting $\tilde T = \chi_{Y \setminus C}T,$ we can write $\lambda [T_C] + [\tilde T]
= 0.$
\sn Let us first assume $dT = 0.$ Then $\tilde T$ defines canonically a current $T'$ on $X$ by
taking the trivial extension of $TÊ\vert X \setminus E.$ By [Sk82], $T'$ is closed. Letting $C_0 \subset E =
\P(N^*_{C \vert Y})$ be a section  with minimal self-intersection, we obtain :
$$ \lambda [T_{C_0}] + [T'] = \mu [l].$$
Obviously $\mu > 0,$ since $X$ is K\"ahler. 
Since $l$ is extremal, we must have $T' = 0,$ i.e. $\tilde T = 0,$ and $\lambda = 0,$
hence $T = 0.$
\sn Now we treat the general case. We define $\tilde T$ and $T'$ as before. We prove that
$\d \dd T' = 0,$ the analogous statement for $\tilde T$ being completely parallel. 
By [Ba94,3.5], we have
$$ \d \dd T' = (\d \dd T \vert X \setminus E)^0 + R,$$
where $( \  )^0$ denotes again the trivial extension and $R$ is a negative current supported on
$E.$ Choosing a negative closed (2,2)-form $u$ on $X$ as in the proof of (3.7) so that $T' + u
\leq 0$ we have furthermore by loc.cit.
$$ \d \dd (T' + u) = (\d \dd (T + u) \vert X \setminus E)^0 + R' = (\d \dd T \vert X
\setminus E)^0 + R',$$
where $R'$ is a positive current supported on $E.$ In total we get $R = R' = 0,$
hence $\d \dd T' = 0.$
\sn It is clear that we still have the equation
$$ \lambdaÊ[C] + [\tilde T] = 0.$$      

\sn Next we define a class $[T'] \in H^{2,2}_{\bf R}(X)$  
ad hoc in the following way.
By duality we define instead a linear form $$[T'] : H^{1,1}_{\bf R}(X) \La {\bf R},$$ 
by setting $[T'](\alpha) = T'(\eta),$ where $\eta$ is any $\d-$ ( $=
\dd-$ closed representative, since the class is real) of $\alpha.$ 
Then there is still an equation
$$ \lambdaÊ[T_{C_0}] Ê+ [T']  = \mu [l] $$
as before. Cupping with $[\omega]$ gives $\mu \geq 0$ and in fact $\mu > 0,$ if
$T \ne 0.$ Noting $[T']Ê\in {\overline {NA}}(X),$ the extremality of $l$ yields 
$\lambda = 0, T' = 0,$ so $T = 0.$    

\bn {\bf Remark} \  Of course one would expect that a stronger version of 3.11 holds, namely that $Y$ is
K\"ahler if and only if $l$ is extremal in ${\overline {NE}}(X).$ 

\bn \bn We go back to our special situation and assume that $Z_0$ is not K\"ahler in our
contraction $\psi : X_0 \La Z_0, $ which blows up the curve $C_0 \subset Z_0.$ Then the 
fibers of $\psi_0$ deform into $X_t$ and we obtain a family $\psi_t : X_t \La Z_t $  for 
small $t.$ Every $\psi_t$ is a blow-down of a divisor $E_t$ which fit together to a divisor
$E.$ In principle we would like to show that by chosing $\varphi_{t_0}$ at the beginning
carefully, $Y_0$ must be K\"ahler. However we do not know how to do this at the moment.
Therefore we examine a simpler situation, namely that $\XX$ has only one extremal ray, so
that there is no choice. In that case the extremal rays should converge to an extremal ray
on $X_0,$ i.e. $Y_0$ should be K\"ahler. Let us still simplify the situation, namely that
all $X_t$ are projective and hence that all $Y_t$ are projective, $t \ne 0.$  Then still
we cannot conclude that $Y_0$ is K\"ahler because the K\"ahler cone on $X_0$ might be
smaller than the K\"ahler cone of the nearby fibers. At least we can state the following
proposition which seems to be of independent interest.

\bn {\bf 3.12 Proposition} \ {\it Let $\pi: \X \La \Delta$ be a family of compact K\"ahler 
threefolds, $E \subset \X$ a family of ruled surfaces over $\Delta$ and $\varphi =
(\varphi_t) : \X \La \Y$ be the simultaneous blow-down of $E$ to the manifold $\Y \La
\Delta.$ Let $h^{1,1}(X_0) = 2.$ Assume that there is a sequence $t_{\nu}$ converging to
$0$ such that the $Y_{t_{\nu}}$ are K\"ahler. Then $Y_0$ is K\"ahler (hence all $Y_{t_{\nu}}$ are
K\"ahler for small $t.$) }

\bn {\bf Proof.} Choose a family $(\omega_t)$ of K\"ahler metrics on $(X_t).$ Assuming
$Y_0$ to be non-K\"ahler and letting $C_t = \varphi_t(E_t),$ we find a positive $\d \dd-$closed
current $\tilde T$ such that 
$$ C_0 + \tilde T = \d{\overline S} + \dd S ,$$
compare (3.11). We define its class  $[\tilde T] \in H^{2,2}_{\bf R}(Y_0)$ as in (3.11).
Then
$[C_0] + [\tilde T] = 0.$ Now consider the current $\varphi_{0*}(\omega_0);$ it defines a class in
$H^{1,1}_{\bf R}(Y_0),$ and we obtain 
$$ ([C_0] + [\tilde T]) \cdot  [\varphi_{0*}(\omega_0)] = 0.$$ 
Let $T'$ be the trivial extension of $\tilde T \vert X_0 \setminus E_0$ to $X_0.$ By the same 
arguments as in (3.11), we have $\d \dd T' = 0$ and define its class 
$[T'] \in H^{2,2}_{\bf R}(X_0).$  Let $l \subset E$ be a fiber of
$\varphi_0.$ Then 
$$ [T'] = (\varphi_0)^*[\tilde T] + \lambda [l],$$
with $\lambda \leq 0,$ by application of (3.7.a,b). Therefore from
$$ [T'] \cdot [\omega_0] = [\tilde T] \cdot [\varphi_{0*}(\omega_0)] + \lambda [l] \cdot [\omega_0] $$
we obtain $[ \tilde T]Ê\cdot [\varphi_{0*}(\omega_0)]Ê\geq 0,$ hence 
$$ [C_0]Ê\cdot [\varphi_{0*}(\omega_0)] \leq 0. \eqno (*)$$ 
Of course the same considerations hold on $Y_t.$ But now $[\varphi_{t_{\nu}}(\omega_{t_{\nu}})]$
is represented by a K\"ahler form, since $h^{1,1}(Y_{t_{\nu}}) = 1.$ Therefore
$$ [C_{t_{\nu}}] \cdot [\varphi_{t_{\nu}}(\omega_{t_{\nu}})] > 0.$$
By continuity we obtain from $(*)$ that
$$ [C_0]Ê\cdot [\varphi_{0*}(\omega_0)] = 0.$$
Since $h^{1,1}(Y_0) = 1,$ this means $[C_0] = 0 $ in $H^4(Y_0,{\bf R}) .$ Since $[C_t]$ is
an integer point, we conclude that $[C_t] = 0$ for all $t,$ contradiction. 

\bn {\bf 3.13 Remark} \ If in (3.12) we have $h^{1,1}(Y_t) \geq 2,$ then we cannot conclude
because we don't necessarily have
$$ [C_t]Ê\cdot [\varphi_{t*}(\omega_t)]Ê> 0.$$
Certainly we can choose some K\"ahler form $\omega'$ on $X_t$ with this property, but this
$\omega'$ might not fit into a family $(\omega_t)$ of K\"ahler metrics on all of $\X.$ 
So it seems not impossible that (3.12) does not hold in general.
\sn Unfortunately (3.12) does not give anything in our special situation, because 
$h^{1,1}(X_0) = 2$ forces $X_0$ to be projective :

\bn {\bf 3.14 Proposition}Ê\ {\it Assume that in (3.12) all the $\XX$ are projective. 
Then all $\YY$ are projective and $X_0$ is projective. }

\bn {\bf Proof.}Ê\ Assume $X_0$ non projective, equivalently $Y_0$ non projective (3.12).
\sn Since $h^{1,1}(\YY) = 1,$ we have $${\rm Pic}(\YY) /{\rm torsion} \simeq {\bf Z}. $$ Fix a big generator
$\O_{\YY}(1)$ of Pic$(\YY)/{\equiv}.$ Then
$$ K_{\YY}Ê\equiv \O_{\YY}(\alpha)$$
for some $\alpha \in {\bf Z}.$ If $\alpha > 0,$ then $\YY$ is of general type and so does
$Y_0.$ Hence $X_0$ is projective. If $\alpha < 0,$ then $\kappa(-K_{\YY}) = 3$ and so does
$\kappa(-K_{\Yy}).$ The conclusion is as before. So we are left with $\alpha = 0.$ 
Then $\chi(\O_{\YY}) = 0$ by Riemann-Roch. $X_0$ being non-algebraic by assumption, we have
$H^2(\O_{X_0}) \ne 0.$ Since $h^3(\O_{\YY}) = h^0(K_{\YY}) \leq 1,$ we conclude that the
irregularity $q(\YY) = h^1(\O_{\YY}) \geq 1.$ Hence we have a non-trivial Albanese map
$$ \alpha : \YY \La {\rm Alb}(\YY).$$
Since $\kappa(\YY) = 0, \alpha $ is surjective [Be83]. This immediately contradicts
$h^{1,1}(\YY) = 1. $ 

\bn But now an examination of the proof of (3.14) shows that it demonstrates at the same
time 

\bn {\bf 3.15 Proposition} \ {\it Assume the situation of (3.14) but instead of an assumption
on $h^{1,1}$ assume $\rho(\XX) = 2$
for some $\nu$ and moreover that $X_0$ is not projective. Then all $Y_t$ are tori, in
particular $Y_0$ is K\"ahler. }

\bn In fact, with the arguments in the proof of (3.14) we arrive at the Albanese map
$\alpha : \YY \La {\rm Alb}(\YY)$ which is onto and has connected fibers. Since $\rho(\YY) = 1,$
dimÊAlb$(\YY) = 3$ and $\alpha $ does not contract any divisor. But since Alb$(\YY)$ is
smooth, $\alpha$ cannot contract any curve either. So $\YY = {\rm Alb}(\YY) $ (alternatively
apply directly the decomposition theorem [Be83]). Hence our claim follows easily.

\bn Some parts of the arguments of this section do not need the condition on algebraic
approximability. Indeed we have  

\bn {\bf 3.16 Theorem} {\it Let $X$ be a compact K\"ahler threefold with $\kappa (X) \geq 0.$
Assume that there
exists a rational curve $C \subset X$ such that $K_X \cdot C < 0.$ Then there exists
a contraction $\varphi : X \La Y$ as described in (3.4)}

\bn {\bf Proof.} As already explained, we construct from the curve $C$ a non-splitting 
family $(C_t)$ of rational curves. 
Then we can apply again the main theorem of Part 1 and are 
through. It was necessary to assume that $\kappa (X) \geq 0$ in order to exclude the
case $-K_X \cdot C_t = 1$  and the $C_t$ fill up a non-normal surface (in which case
we only have informations if $X$ can be approximated algebraically or if $\kappa(X) = 
- \infty).$     

\bn \bn \bn 
{\bf 4. An Openness Theorem  } 
\bn In this section we shall prove
\bn 
{\bf 4.1 Theorem} \ {\it Let $\pi: \X \La \Delta$ be an algebraic approximation of the 
non-algebraic compact K\"ahler threefold $X_0.$ Let $(t_{\nu}) \to 0$ be a sequence such
that $\XX$ is projective for every $\nu.$ If $K_{X_0}$ is not nef, then all $K_{\XX}$
are not nef for $\nu \gg 0.$}

\bn We will also prove that on threefolds with $\kappa (X) = 2$ and $K_X$ not
nef, there exists a rational curve $C$ with $K_X \cdot C < 0$ (theorem 4.10 ). If $\kappa (X) = 1,$
then we can at least prove that there is {\it some} curve $C$ with $K_X \cdot C < 0$
(theorem 4.15).   

\bn {\bf 4.2 Remark} \ If $\X = (X_t)$ is a family of smooth projective threefolds such
that $K_{X_0}$ is not nef, then all $K_{X_t}$ are not nef. This follows e.g. by deforming
extremal rational curves in $X_0$ to $X_t.$ For details and the higherdimensional context
see [AP95].

\bn The proof of (4.1) will be given in several steps according to the possible values of
the Kodaira dimension. We fix a sequence $(t_{\nu})$ as above and assume that $K_{\XX}$ is
nef for all $\nu.$ 

\bn {\bf 4.3 Proposition}Ê\  {\it We have $\kappa(X_t) \geq 0$ for all $t \in \Delta.$ }

\bn {\bf Proof.} \ Assume that $\kappa(X_{t_0}) = -\infty$ for some $t_0 \in \Delta.$
If $X_{t_0}$ is projective, then $X_{t_0}$ is uniruled [Mo88],[Mi88],[Ka92], hence all
$X_t$ are uniruled, contradiction. 
\sn So let $X_{t_0}$ be non-projective; again we want to show that $X_{t_0}$ is uniruled
which ends the proof (this proof works also in the projective case, making the argument 
independent from the deep papers cited above). Since $K_{\XX}$ is nef, we have
$$ \chi(\XX,\O_{\XX}) \leq 0 $$
by [Mi87]. Therefore $\chi(\O_{X_{t_0}}) \leq 0.$ Since $h^3(\O_{t_0}) = h^0(K_{X_{t_0}}) = 0$ by our assumption, we conclude
that the irregularity $$q(X_{t_0}) = h^1(\O_{X_{t_0}}) \geq 1.$$ So we have a non-trivial Albanese
map
$$ \alpha : X_{t_0} \La {\rm Alb}(X_{t_0}).$$ If dim $\alpha(X_{t_0}) = 1,$ then by
$C_{3,1}$ (see [Fj78]) the general fiber $F$ has  $\kappa(F) = - \infty,$ hence $X_{t_0}$
is uniruled. If dim $\alpha(X_{t_0}) = 2,$ then by $C_{3,2}$ [Ka81], the general fiber of $\alpha$
is ${\bf P}_1$ and again $X_{t_0}$ is uniruled. If finally dim $\alpha (X_{t_0}) = 3,$ then
$\kappa(X_{t_0}) \geq 0, $ a contradiction. 

\bn Next we deal with two easy cases. Observe first that $\kappa(\Xx)$ does not depend on
$\nu_0$ since by [Ka92]
$$ \kappa(\Xx) = {\rm max} \ \{m \in {\bf N}Ê\cup \{0\} \vert K^m_{\XX}ÊÊ\not \equiv 0 \}.$$

\bn {\bf 4.4 Proposition} \ {\it (1) If $K_{\XX} \equiv 0$ for some $\nu$, then $K_{X_0} \equiv 0$
and $X_0$ is  up to finite etale cover either a torus or a product of an elliptic curve with a
K3 surface.
\sn (2) If $K^3_{\XX} > 0$ (i.e. $K_{\XX}$ is big) for some $\nu,$ then $\kappa(X_0) = 3$ 
and $X_0$ is projective. }

\bn {\bf Proof.} \ (1) Just notice that $c_1(\XX) = 0$ in $H^2(\XX,{\bf Q})$ implies that
$c_1(X_0) = 0$ in $H^2(X_0,{\bf Q}).$ Hence we conclude by [Be83].
\sn (2) By Grauert-Riemenschneider or Kawamata-Viehweg vanishing, we have for $q \geq 1$ that
$$ H^q(\XX,mK_{\XX}) = 0 \eqno (*)$$
for all $\nu$ and all $m \geq 2.$ By the coherence of $R^q\pi_*(mK_{\X}),$ equation (*) holds for
all $t \in \Delta \setminus S,$ with $S \subset \Delta $ a proper analytic subset. Passing to a local
1-dimensional submanifold of $\Delta$ through $0$, we may assume that (*) holds for all $t \ne 0.$ 
By Riemann-Roch we conclude $$ \chi(X_t,mK_{X_t}) = h^0(X_t,mK_{X_t}) \sim m^3$$ for $t \ne 0,$
and this function in $t$ is constant. By semi-continuity of cohomology it follows
$$ h^0(X_0,mK_{X_0}) \sim m^3,$$
so that $X_0$ is Moishezon. Since $X_0$ is K\"ahler, it must be projective, contradiction.

\bn We are therefore reduced to the intermediate cases $\kappa(X_{t_{\nu}}) = 1$ or $2.$ These
cases are much more complicated and we shall need some preparations. The first lemma is
probably well-known but apparently never stated explicitly.

\bn {\bf 4.5 Lemma} \ {\it  Let $X$ be a compact K\"ahler manifold and $L$ a line bundle on
$X.$ Then $L$ is nef if and only if $c_1(L)$ is in the closure of the K\"ahler cone.}

\bn {\bf Proof.} \ (a) First assume that $L$ is nef. Take an arbitrary K\"ahler metric 
$\omega $ on $X.$  Then we must show that $c_1(L) + [\omega]$ is represented by a K\"ahler
metric. For $0 < \varepsilon < 1 $ choose a metric $h_{\varepsilon}$ on $L$ with curvature
$$\Theta_{L,h_{\varepsilon}} \geq - \epsilon \omega.$$
Then $\Theta_{L,h_{\varepsilon}} + \omega \geq (1-\varepsilon) \omega,$ hence 
$\Theta_{L,h_{\varepsilon}} + \omega$ is a K\"ahler form representing $c_1(L) + [\omega]$.
\sn (b) Now suppose that $c_1(L)$ is in the closure of the K\"ahler cone. Fix a K\"ahler
metric $\omega$ on $X.$ Let $\varepsilon > 0.$ By [Su88] we find a K\"ahler metric 
$\omega_{\varepsilon}$ such that 
$$ [\omega_{\varepsilon} - \varepsilon \omega] = c_1(L) . $$
Choose a metric $h_{\varepsilon}$ on $L$ with curvature $\Theta_{L,h_{\varepsilon}} =
\omega_{\varepsilon} - \varepsilon \omega.$ Then $\Theta_{L,h_{\varepsilon}}Ê\geq
- \varepsilon \omega.$ Since $\varepsilon > 0$ was arbitrary, $L$ is nef.  

\bn {\bf 4.6 Lemma} {\it Let $X$ be a smooth compact K\"ahler surface or threefold and $A \subset X$ a compact
submanifold. Let $\pi: \hat X \La X$ be the blow-up of $A.$ Let $L$ be a line bundle on $X.$
Then $L$ is nef if and only  if $\pi^*(L)$ is nef. }

\bn {\bf Proof.} For simplicity of notations we only treat the threefold case.
One direction being obvious, we assume that $\pi^*(L)$ is nef. Let $KC(X)$
and $KC(\hat X)$ denotes the K\"ahler cones on $X$ resp. $\hat X.$ By the easy part of (4.5)
we have 
$$ c_1(\pi^*(L)) \in KC(\hat X).$$
We claim that this implies 
$$ c_1(L) \in KC(X). \eqno (1)$$
Then the proof is finished by applying the difficult part of (4.5). 
\sn For the proof of (1) we need to show that
$$ T \cdot c_1(L) \geq 0 \eqno (2)$$
for every positive closed current of bidegree $(2,2).$ 
We know that $\hat T \cdot c_1(\pi^*(L)) \geq 0$ for all closed positive currents
$\hat T$ on $\hat X$ of bidegree $(2,2).$ Therefore the proof of (2) is reduced to
the following. 
\sn (3) For every closed positive current $T$ on $X$ of bidegree $(2,2)$ there exists a
closed positive current $\hat T$ on $\hat X$ such that
$$ \pi_*(\hat T) = T.$$ 
(Note that the dual cone to the K\"ahler cone is the cone of positive closed (1,1)-currents).
\sn Write
$$  T = \chi_AT + \chi_{X \setminus A}T.$$
If ${\rm dim}A = 0,$ then $\chi_AT = 0,$ otherwise $\chi_AT = \lambda T_A$ with 
$\lambda \geq 0.$
\sn Now let $T' = T \vert X \setminus A.$ Identifying $X \setminus A$ with $\hat X \setminus
\hat A,$ we  can consider the trivial extension $T_1$ of $T'$ on $\hat X.$ By [Sk82] this
is a closed positive current and obviously $$\pi_*(T_1) = \chi_{X \setminus A}T.$$ We are thus
reduced to dealing with $\chi_AT$ and may assume ${\rm dim}A = 1.$ But then we take a
section $$C \subset \hat A = \P(N^*_{A \vert X})$$ and let $T_2 = \lambda T_C.$ Hence
$\pi_*(T_2) = \chi_AT$ and we are done, putting $\hat T = T_1 + T_2.$ 

\bn {\bf Remark.} Of course (4.6) should be true in any dimension.  

\bn Lemma 4.6 enables us to define nefness on any reduced compact complex space of dimension at
most three

\bn {\bf 4.7 Definition} ÊLet $X$ be a reduced compact complex space, ${\rm dim}ÊX \leq 3.$ 
Let $L$ be a
line bundle on $X.$ Then $L$ is {\it nef}, if there exists a desingularisation $\pi: \hat X \La X$
such that $\pi^*(L)$ is nef.

\bn By (4.6) the definition does not depend on the choice of the resolution.

\bn {\bf 4.8 Problem} \ There is another natural way to define nefness for a line bundle
$L$ on a reduced compact complex space $X:$ we require that for any $\varepsilon > 0$ there
exists a metric $h_{\varepsilon}$ on $L$ such that the curvature satisfies
$\Theta_{L,h_{\varepsilon}}Ê\geq - \varepsilon \omega $ on the smooth part of $X,$
where $\omega$ is a fixed positive (1,1)-form on $X.$ 
\sn These both notions of nefness should coincide. 

\bn The next lemma will be very important for the main results of this section; it is of course
trivial in the projective case.  

\bn {\bf 4.9 Lemma } \ {\it Let $X$ be a compact K\"ahler manifold with ${\rm dim} X \leq 4$ and let $L$ be a line bundle on $X.$
Assume that there is an effective divisor $D$ on $X$ such that $L =
\O_X(D).$ If $L \vert D$ is nef, then $L$ is nef.}

.

\bn {\bf Proof.}Ê\ Note that the K\"ahler cone is the dual cone of the cone
of of classes of positive closed currents of bidegree $(n-1,n-1) $ on $X.$ 
So by (4.5) we only need to show
$$ \int_X c_1(L) \wedge T \geq 0 \eqno(*) $$
for every positive closed $(n-1,n-1)-$ current $T.$ Of course we may assume that
$D$ is a reduced prime divisor.
By Skoda [Sk82], the current
$\chi_DT$ is closed, hence we obtain from (*) :
$$\int_X c_1(L) \wedge T = \int_X c_1(L) \wedge \chi_DT + \int_X c_1(L) \wedge \chi_{X
\setminus D}ÊT. \eqno(**) $$
Appyling (4.6) to an embedded resolution of singularities for $D$, we may assume $D$ to be
smooth from the beginning. Then it is an easy exercise
to conclude 
$$ \int_X c_1(L) \wedge \chi_DT \geq 0 $$  
(extend metrics on $L \vert D$ to a neighborhood). 
\sn Now the second term in (**) is also non-negative by (3.7.a), hence (*) follows and the
claim is proved.

\bn The assumption ${\rm dim} X \leq 4$ was used only for being able to apply (4.6).

\bn {\bf 4.10 Theorem}  {\it Let $X$ be a compact K\"ahler threefold with $\kappa (X) = 2.$
If $K_X$ is not nef, there exists a rational curve $C$ with $K_X \cdot C < 0.$}

\bn {\bf Proof.} Let $m \gg 0$ so that $V = \vert mK_X \vert $ defines a rational map
$$ f :  X \rightharpoonup Y$$ to a projective surface $Y.$ We choose a sequence of blow-ups
$\pi : \hat X \La X$ such that the induced map $\hat f : \hat X \La Y$ is holomorphic.   
We note that $f$ is almost holomorphic, i.e. there is a non-empty Zariski open set $U \subset X$
such that $f \vert U$ is holomorphic and proper. In fact, otherwise the exceptional set of
$\pi$ would contain a component lying surjectively over $S$ and therefore $X$ would be
projective by [Ca81].   
\sn Let $D_0 \in \vert mK_X \vert $ be
a general element.  By (4.9) $K_{X} \vert D_0$ is not nef. 
Observe that 
$V$ must have a positive dimensional base locus, otherwise $K_{X}$ would be nef. 
We claim that $V$ must even have fixed components:
\sn otherwise the general $D_0$ is irreducible. 
It follows from (4.14) below (applied to a desingularisation) that there is a curve $C \subset D_0$ with $K_{X} \cdot C
< 0$. This is impossible since then $C$ would deform
in $X_0$ in a 1-dimensional family [Ko91] to fill up $D_0$ which is clearly absurd.

\sn Hence $V$ has fixed components, $D_0$ is reducible and we can write
$$ D_0 = B + \sum_{i}Ê\lambda_i A_i$$
where $B$ is the movable part of $D_0,$ hence irreducible and with multiplicity $1$, and the
$A_i$ are the fixed components. 

\bn {\bf (4.10.1)} Our first aim is to show that there must be an irreducible curve $C \subset X$
(not necessarily rational) with $K_{X}Ê\cdot C < 0.$ By the arguments proving the existence of the fixed components it
follows also that we can find some $i_0$ such that $K_{X} \vert A_{i_0}$ is not nef (apply
the above argument for $B$ instead of $D_0).$ Moreover
it suffices in order to get a contradiction to show that $A_{i_0}$ has a non-constant meromorphic
function and that the case (4.13(2)) does not occur. Let $A = A_{i_0}$ for simplicity of
notations, and let $\mu = \lambda_{i_0}.$ 
We can choose $ \pi: \hat X \La X$ 
such that the strict transform $\hat A$ of $A$ is smooth. We investigate
the structure of $\hat A$ by distinguishing cases according to dim$\hat f(\hat A).$ 
\sn (a) dim$\hat f(\hat A) = 0.$
\sn Let $\hat F$ denote the fiber of $\hat f$ containing $\hat A,$ equipped with the natural
structure. Write $\hat F = \lambda {\hat A} + R.$ Then the conormal bundle $N^*_{\hat F \vert 
\hat X}Ê\vert \hat A$ is generated by global sections, but is clearly not trivial, hence
$$ \kappa(N^*_{\hat FÊ\vert \hat X}Ê\vert \hat A) \geq 1,$$
in particular $a(A) \geq 1.$ Therefore we only need to exclude the case (4.13(2)) to conclude.
So let $\hat h : \hat A \La \hat C$ be the algebraic reduction of $\hat A$ and assume 
$$ \pi^*(K_{X}Ê\vert A) \equiv \hat h^*(\hat G)$$
with $\hat G^*$ nef and not numerically trivial, hence ample. Clearly
$\hat h$ descends to a map $h : A \La C$ with
$g : \hat C \La C$ being the normalisation. Then $K_{X}Ê\vert A \equiv  h^*(G)$ with $G^*$ ample
on $C.$ Write 
$$mK_{X} = D_0 = \mu A + M.$$
The adjunction formula gives immediately the follwoing general fact. If $X$ is a complex manifold,
$A,B \subset X$ hypersurfaces in $X$ and $D = A + B,$ then $K_A = K_D \vert A - B \vert A.$ 
Thus we have in our situation
$$ K_{\mu A} \vert A = mK_{X}Ê\vert A - M',$$
where $M'$ is effective on $A$ and supported on $M \cap A.$
On the other hand by ordinary adjunction
$$ K_{\mu A}Ê\vert A = K_{X}Ê\vert A + N_{\mu AÊ\vert X}Ê\vert A,$$
hence in total $N^*_{\mu A \vert X}Ê\vert A = N^{*\mu}_A$ is effective by the ampleness of $G$ and 
effectivity of $M'.$ Therefore from $K_A = K_{X} \vert AÊ+ N_A$ we see that $\kappa(K_A) = 
- \infty.$ Since $K_{\hat A}$ is a subsheaf of $\pi^*(K_A)$, it follows $\kappa(\hat A) 
= - \infty$ and $\hat A$ is algebraic, contradiction. 

\bn (b) dim$\hat f(\hat A) = 1.$
\sn Now it is obvious that $a(\hat A) \geq 1.$ The remaining case (4.13(2)) is excluded as
in (a).

\sn (c) dim$\hat f(\hat A) = 2.$
\sn Then $\hat f\vert \hat A$ is onto $Y$. However $Y$ is algebraic and so does $A$. Then
our claim is obvious.
\sn (4.10.1) is thus completely proved. 

\bn {\bf (4.10.2)} We next claim that there is a {\bf rational} curve $C \subset X$ such
that $K_{X} \cdot C < 0.$ 
\sn If such a rational $C$ does not exist then we find an irrational curve $C \subset X$
with $K_{X} \cdot C < 0,$ hence $C$ deforms (with irreducible parameter space) in an 
at least 1-dimensional family $(C_t)_{t \in T}$
[Ko91] with $C = C_0.$ Moreover we may assume that no deformation of $C$ splits. We choose
$T$ maximal (but of course irreducible, as always). Then $\bigcup_{t \in T} C_t $ is a fixed
component of $V,$ say $A.$ Write a general element of $\vert mK_X \vert$ as  
$$ D_0 = m_1 A + R$$
so that $R$ does not contain $A.$
Then by adjunction
$$ m_1 K_A = m_1 K_{X} \vert A + m_1 A \vert A \subset m_1 K_{X}Ê\vert A + D_0 \vert A 
= (m_1 + m)K_{X}Ê\vert A. $$
Therefore we obtain the following basic inequality
$$ K_A \cdot B \leq (1 + {Êm \over m_1})ÊK_{X}Ê\cdot B < K_{X}Ê\cdot B \eqno (I) $$
for all but finitely many curves $B \subset A$ with $K_X \cdot B < 0.$  Hence
$$ K_A \cdot C_t \leq -2.  \eqno (*)$$
Let $\nu : \tilde A \La A$ be the normalisation of $A$ and $\pi' : \tilde {\tilde A} \La \tilde A$ the minimal
desingularisation. Then we have
$$ K_{ \tilde A} = \nu^*(K_A) - E$$
with $E$ an effective Weil divisor supported exactly on the preimage of the non-normal locus
of $A,$ cp. [Mo82]. Moreover
$$ K_{\tilde {\tilde  A}}Ê= \pi'^*(K_{\tilde A}) - E'$$
with an effective divisor $E'$ which is $0$ if and only if $\tilde A$ has only rational double
points. 
\sn Let $t \in T$ be general and $\tilde C_t$ be the strict transform of $C_t$ in $\tilde {\tilde A}.$ 
Then (*) yields
$$ K_{\tilde {\tilde  A}}Ê\cdot {\tilde C_t} \leq -2. \eqno (**)$$
\sn The general fiber $\hat F$ of $\hat f$ (which is an elliptic curve) induces a unique maximal 
family $(\hat F_t)_{t \in  \hat T}$ with graph
$$q: \C \La \hat T, p: \C \La \hat X$$ (taking closure in the cycle space). The map $p$ is clearly
bimeromorphic, $f$ being almost holomorphic, and $q$ is an elliptic fibration. Of course we may
assume $\hat T$ smooth. Let $\overline A$ denote the strict transform of $\tilde {\tilde A}$ in $\C;$ we
may assume $\overline A$ smooth (by possibly choosing an embedded resolution of singularities
and then flattening $q$). Let $g = q \vert \overline A : \overline A \La D \subset \hat T.$ By (**)
we have $$\kappa (\tilde {\tilde A}) = - \infty.$$ Note also that $D$ is a
rational or elliptic curve. 
\sn (a) Let us first assume that the genus $g(C_t) \geq 2$ ($=$ genus of the normalisation).
Our aim is to construct a new family of elliptic or rational curves $(C'_t)$ such that
$K_X \cdot C'_t < 0.$
 
\sn Let $\hat G_t = \hat F_t \cap \overline A$ for $t \in D.$ Since $f$ is almost holomorphic,
we have
$$ p^*\pi^*K_X \cdot \hat F_t = 0.$$
Now take a component $\hat C \subset \hat F_t$ with ${\rm dim} \  C = {\rm dim} \  \pi \circ p (\hat C)
= 1.$ If 
$$ p^* \pi^* K_X \cdot \hat C > 0,$$
then we find another $\hat C' \subset \hat F_t$ with ${\rm dim} \ Ê\pi \circ p (\hat C') = 1$
such that $$p^* \pi^*(K_X) \cdot  \hat C' < 0.$$ Let $C' = \pi \circ p (\hat C').$ Then
$K_X \cdot C' < 0$ and we have found a family of elliptic or rational curves $(C'_t)$
with $K_X \cdot C'_t < 0.$ \sn Therefore we may assume $K_X \cdot C = 0$ for any choice of
$\hat C \subset \hat F_t, t \in D$ and $C = \pi \circ p(\hat C).$ We shall assume now that the 
general $G_t$ is elliptic,
the case that $G_t$ is rational being even easier (and therefore omitted). We consider the possibly not relatively
minimal fibration $$g: \overline A \La D \simeq \P_1.$$
The curve $D$ is rational since $\kappa (\overline A) = - \infty.$  
Furthermore we have a "ruling" $h: \overline A \La E,$ so that the general fiber of $h$ is
a smooth rational curve and every fiber is a tree of $\P_1.$ Let $$L = p^*  \pi^*(K_X).$$ 
Then  $L \cdot  \hat C = 0$ for every component of every fiber of $g.$ Since $f$ is almost
holomorphic, we even have $$ p^*\pi^*(K_X) \cdot \hat F_t = \O_{\hat F_t}$$ for general $t$. From
semi-continuity it follows then easily that $L \vert \hat C = \O_{\hat C}$ for all $C,$ possibly up to
a torsion line bundle (in case $C$ appears with multiplicity in the fiber). After eventually passing
from $L$ to $L^m$ we get  
$$ L \equiv g^*(\O_D(a))$$ 
for some integer $a.$ 
Now let $C_t$ be a general element of our family $(C_t)$ and let $\overline C_t$ be its strict
transform in $\overline A.$ Then 
$$ L \cdot \overline C_t < 0.$$ 
Since $\overline C_t$ is a multi-section of $g,$ we conclude $a < 0.$ Now consider a general
fiber $l$ of $h.$ Then we conclude $L \cdot l < 0,$ and therefore
$$ K_X \cdot \pi p(l) < 0.$$ Therefore the images of $l$ define the family $C'_t$ of rational
curves. 
\sn (b) We are now reduced to the case that $(C_t)_{t \in T}$ is a family of elliptic curves and moreover
we may assume that the family is non-splitting in the sense that no component of any
member is a rational curve (but multiples of elliptic curves are allowed). We take $T$ of course
maximal and assume $T$ normal.  Denote by $\D$ the normalisation of the graph of the family with projections
$p_1: \D \La A$ and $q_1: \D \La T.$ \sn First assume ${\rm dim} T \geq 2.$ Note for the following that
$T$ is algebraic ( = Moishezon) since $A$ is algebraic and $A = \bigcup C_t.$ Choose a general point $x \in A$
and introduce
$$ T(x) = \{ t \in T \vert x  \in C_t \}.$$ 
Then ${\rm dim} T(x) \geq 1.$ Choose an irreducible curve $\Delta \subset T(x).$ 
Let $p_{\Delta}Ê: \D_{\Delta}Ê\La A$ be the normalisation of the induced graph with
projection $q_{\Delta}: \D_{\Delta}Ê\La \Delta.$ Then $q_{\Delta}$ admits a multi-section
$Z$ with ${\rm dim} \  p_{\Delta}(Z) = 0,$ so that $Z \subset \D_{\Delta}$ is exceptional.
This contradicts Sublemma (4.10.a). 
\sn We are therefore left with the case ${\rm dim} T = 1.$ Since $K_A$ is a subsheaf of
$K_X \vert A,$ the inequality
$$ K_X \cdot C_t \leq -1$$
implies
$$ K_A \cdot C_t \leq -2$$ 
by virtue of the inequality (*) proven above.
\sn
For general $t,$ let $\hat C_t$ be the strict transform of $C_t$ in $\hat A,$ as before.
Then $$ K_{\hat A}Ê\cdot \hat C_t \leq -2.$$
Hence $$h^0(N_{\hat C_t}) \geq 2$$ 
whereas 
$$H^1(N_{\hat C_t}) = 0.$$ 
Consequently $\hat C_t$ moves in an at least 2-dimensional family, contradiction. 
\vfill \eject
 
\sn The proof of (4.10) is now complete modulo the following

\bn {\bf 4.10.a Sublemma} {\it Let $S$ be a smooth projective surface and $f: S \La C$
a holomorphic elliptic fibration onto a smooth curve $C.$ Assume that the only singular fibers are
multiples of elliptic curves.  Then there is no irreducible
curve $A \subset S$ with $A^2 < 0.$ }

\bn {\bf Proof.} Let $g$ be the genus of $C.$ By adjunction 
$$ {\rm deg} K_A = K_S \cdot A + A^2.$$
The canonical bundle formula for elliptic fibrations (see e.g. [BPV84,chap.5, 12.1
and chap.3, 18.2]) gives $K_S \equiv f^*(K_C).$ Let $d = {\rm deg} f \vert A.$
Then $$K_S \cdot A = d(2g-2).$$
Let $\nu : \tilde A \La A$ be the normalisation of $A;$ then 
${\rm deg}ÊK_{\tilde A} \leq  {\rm deg} K_A.$ Now the formula of Riemann-Hurwitz 
yields 
$$ {\rm deg}ÊK_{\tilde A}Ê= d(2g-2) + {\rm deg}ÊR,$$
where $R$ is the ramification divisor of $ \tilde A \La C.$ 
Putting things together we obtain 
$$ {\rm deg}ÊK_A = d(2g-2) + A^2 \geq  d(2g-2) + {\rm deg}ÊR,$$
hence ${\rm deg} RÊ\leq A^2 < 0,$ which is absurd.

\bn {\bf 4.11 Corollary} {\it Let $X$ be a compact K\"ahler threefold with $\kappa (X) = 2.$
Assume that $K_X$ is not nef. Then there exists a birational extremal contraction  $f : X \La Y$
(i.e. of type (c) or (d) in (3.4)).}

\bn {\bf Proof.} 
By (4.10) exists a rational curve $C \subset X_0$ with $K_{X_0}Ê\cdot C < 0;$ then 
apply (3.16).

\bn {\bf (4.12)} In particular we have proved (4.1) in case $K_{X_0}^3 = 0, K_{X_0}^2 \ne 0.$
In fact, if $ K_{X_0}$ would not be nef, then we find a rational curve $C \subset X_0$
with $K_{X_0}Ê\cdot C < 0$ and which can be deformed to $X_t.$ In fact, taking over the notations
of the last proof, we conclude that in all cases but the case of ${\bf P}_2$ with $N_S = \O(-2)$ we conclude that 
$$ N_{C_t \vert \X} \simeq \bigoplus \O(a_i) $$
with $a_i \geq -1,$ so that the deformations of $C_t$ in $\X$ are unobstructed. Moreover
$$ h^0(\X,N_{C_t \vert \X}) > h^0(X_0,N_{C_t \vert X_0}), $$
so that the $ C_t$ can be deformed to the neighbouring fibers $X_s,$ so that $K_{X_s}$ cannot
be nef. In the remaining case we see by a similar argument that we can deform directly the
${\bf P}_2$ out of $X_0$ and are also done.

\bn In the proof (4.10) we made use of the following

\bn {\bf 4.13 Proposition}
{\it Let $S$ be a smooth compact K\"ahler surface with algebraic dimension $a(S) = 1.$
Let $f: S \La A$ be the algebraic reduction to the smooth curve $A.$ Let $L$ be a line bundle
on $S.$ Let $\equiv$ denote topological equivalence. Then
\item {(1)}Ê$L$ is nef if and only if $L \equiv f^*(G) $ with $G$ nef
{\item {(2)} $L$ is not nef if and only if either
         \itemitem {(a)} there is a curve $C \subset S$  with  $\L \cdot C < 0$
          \itemitem {(b)} $L \equiv f^*(G'), G'^*$ ample. }}

\bn {\bf Proof.} (1) One direction being obvious, we let $L$ be nef. First we show that
$L \cdot C = 0$ for every curve $C \subset S.$ In fact, assume $L \cdot C > 0$ for some $C.$
Then
$$ (kL+C)^2 = k^2L^2 + 2kL \cdot C + C^2.$$
Since $L^2 \geq 0,$ we obtain by [DPS94] $ (kL+C)^2 > 0$ for $k \gg 0,$ hence $S$ would be algebraic.
\sn Note that if $L  \equiv f^*(G),$ it is clear that $G$
must be nef (otherwise we obtain a nef line bundle whose dual has a section with zeroes
which is impossible by [DPS94]). \sn So assume now $L \cdot C = 0$ for all curves $C.$ By
(4.6) we may assume $S$ minimal.
\sn (a) First let $\kappa (S) = 1.$ Since $L \cdot K_S = 0,$ Riemann-Roch gives
$$ \chi(S,mL) = \chi(S,\O_S) > 0.$$
Hence $h^0(mL) > 0$ or $h^0(mL^* + K_S) > 0.$ In the first case we conclude that
$mL \vert F = \O_F$ for every fiber of the algebraic reduction $f,$ hence $mL = f^*(G).$
In the second case we can write
$$K_S = mL + D$$
with $D$ effective. Since $D \cdot F = 0,$ we get $D \vert F = \O_F,$ hence
$$ \lambda L = f^*(G).$$
\sn (b) Now let $\kappa (S) = 0.$ Then either $S$ is a K3 surface or $S$ is a torus.
If $S$ is a K3 surface, then $\chi(S,\O_S) = 2$ and we can conclude as in (a). So
let $S$ be a torus. Then $f$ is an elliptic fiber bundle over an elliptic curve. The line
bundle $L$ defines a section 
$$s \in H^0(C,R^1f_*(\O_S^*)).$$
Fix a point $y_0 \in C.$ Then we find a topologically trivial line bundle $H$ on $S$
such that $$L \vert f^{-1}(y_0) \otimes H \vert f^{-1}(y_0) \simeq \O.$$
Since $S$ is not algebraic, we automatically have
$$ L \vert f^{-1}(y) \otimes H \vert f^{-1}(y)Ê\simeq \O $$
for all $y \in C;$ otherwise we would get a multisection of $f.$
This implies our claim.

\sn (2) Again one direction is clear. So let $L$ be not nef and assume $L \not \equiv f^*(G').$ 
Then by the proof of (1), there is a curve $C_0$ with $L \cdot C_0 \ne 0.$ If our claim would
be false, then $L \cdot C_0 > 0$ and $L \cdot C \geq 0$ for every curve $C \subset S.$ 
Let $F$ be the general fiber of $f.$ Then 
$$ L \cdot F = 0,$$
otherwise $(L+kF)^2 > 0$ for large $k.$ Since dim$f(C_0) = 0,$ we therefore find 
$C_1 \subset f^{-1}f(C_0)$ with $L \cdot C_1 < 0,$ contradiction.

\bn {\bf 4.14 Corollary} {\it Let $X$ be a non-algebraic compact K\"ahler surface. Let $L$
be a line bundle on $X$ such that $L \cdot C \geq 0$ for all curves $C \subset X.$ 
Assume moreover that some power $L^m$ has a section. Then $L$ is nef. }

\bn {\bf Proof.} (a) First we assume that $a(X) =1$ and let $f: X \La C$ be the algebraic
reduction which is an elliptic fibration. Assume that $L$ is not nef. By (4.13) we find
(posssibly after passing to a multiple of $L$) a topologically trivial line bundle $G$ on $X$ such that
$$ L = f^*(H) \otimes G $$
with some negative line bundle $H$ on $C.$ We may assume that already $L$ has a section.
We conclude that
$$ H^0(C,H \otimes f_*(G)) \ne 0. \eqno (*)$$
Note that $f_*(G)$ is a torsion free sheaf, hence locally free. Since $G$ is topologically
trivial, we have $f_*(G) \ne 0$ if and only if $G \vert F = \O_F$ for the general fiber
$F$ of $f.$ On the other hand $f_*(G)$ must be non-zero by (*). Moreover (*) proves
that $f_*(G)$ and hence $G$ has a section (we may a priori assume that $H^*$ is effective).
Now this section cannot have zeroes and therefore $G = \O_X.$ But then (*) gets absurd.

\sn (b) Assume now that $a(X) = 0.$ Let $\pi : X \La X'$ be the map to the minimal model.
We can write
$$ L = \sum \lambda_i E_i + \pi ^*(L')$$
where the $E_i$ are the exceptional components of $\pi$ and $\lambda_i $ are integers.
Since $L$ is algebraically nef, it is clear that all $\lambda_i \leq 0.$ Therefore $L'$
is immediately seen to be algebraically nef. If $X'$ is a torus, then anyway $L'$ is
topologically trivial and we conclude, e.g. by effectivity of $L^m,$ that all $\lambda_i = 0.$
If $X'$ is a K3 surface, then $H^0(L'^m) \ne 0,$ hence $L'^m = \sum a_i C_i$ with
$C_i$ being $(-2)-$curves (the only curves in a K3 surface without meromorphic functions).
But then $L'$ is not algebraically nef unless all $a_i = 0.$ As in the torus case we
now conclude $\lambda_i = 0$ and $L$ is nef.

\bn {\bf 4.15 Theorem}  {\it Let $X$ be a smooth compact K\"ahler threefold with $\kappa (X) = 1.$
Assume $K_X$ is not nef. Then there exists an irreducible curve $C \subset X$ such that $K_X \cdot C
< 0.$ }

\bn {\bf Proof.}   
We follow the lines of (4.10.1). We assume that $K_X$ is algebraically nef and show
that $K_X $ is nef. Let $$f: X \rightharpoonup C$$
be the meromorphic map defined by $\vert mK_X \vert.$ Let $\pi: \hat X \La X$ be a sequence
of blow-ups such that the induced map $\hat f : \hat X \La C$ is a morphism. Write again
$$ mK_X = B + \sum \lambda_i A_i,$$
where $B$ is the movable part. Then $mK_X \vert B$ is effective. Hence by (4.14) we conclude
that $K_X \vert B$ is nef. It only remains to show that $K_X \vert A_i$ is nef for all $i.$
Fix some $i$ and let $A = A_i.$ If $a(A) \geq 1,$ the arguments of (4.10.1) work. So
assume $a(A) = 0.$ Let $\hat A$ be the strict transform of $A$ in $\hat X;$ we may assume
$\hat A$ smooth. Let $\hat F$
denote the fiber of $\hat f$ containing $\hat A$ and write
$$ \hat F = \lambda \hat A + R.$$
It follows from the exact sequence
$$ 0 \La \O \oplus \O = N^*_{\hat F \vert \hat X}Ê\vert \hat A \La N^*_{\lambda \hat A \vert \hat X}
\vert \hat A \La N^*_{\lambda \hat AÊ\vert \hat F}Ê\vert \hat A \La 0 $$
that $$H^0(\hat A, {\rm det}N^{*\lambda}_{\hat A}) \ne 0.$$
Since $\kappa (\hat A) = 0$ ($\hat A$ being bimeromorphic to a torus or a K3-surface),
it follows from the adjunction formula that $$H^0(\hat A, \lambda K_{\hat X}Ê\vert \hat A) \ne 0.$$
Let $g: \tilde A \La A$ be the normalisation of $A$, then $\pi \vert \hat A = g \circ h$
with $h: \hat A \La \tilde A$ the induced map. Hence
$$ H^0(\hat A, g^*h^*(\lambda K_X) - \lambda \sum \mu_i C_i) \ne 0, $$
where $$K_{\hat X}Ê= \pi^* (K_X) + \sum \mu _i E_i$$
and $C_i = \hat A \cap E_i.$ Then we obtain by applying $h_*$ that 
$$ H^0(\tilde A, g^*(\lambda K_X)) \ne 0.$$
Now apply (4.14)!

\bn {\bf 4.16 Theorem} {\it Let $X$ be a compact K\"ahler threefold with $\kappa (X)  = 1.$
Assume that $K_X$ is not nef and $X$ is algebraically approximable. Then there exists a
rational curve $C \subset X$ with $K_X \cdot C < 0.$} 

\bn {\bf Proof.} If $K_{\XX}$ is not nef for some $\nu,$ then the assertion  is clear (sect. 3). So assume
$K_{\XX}$ nef. We have already seen that this implies $K_X^2 = 0.$ This is the only conclusion
we draw from the algebraic approximability. Write as in (4.15)
$$ mK_X = B + \sum_{i = 1}^k \lambda_i A_i.$$ 
Let $\L = \O_X(B)$ and let $f: X \rightharpoonup S$ be the meromorphic map defined
by $H^0(X,\L)$ to the curve $S.$ 
\sn (1) First let us assume that $f$ is holomorphic, i.e. $H^0(X,\L)$ is base point free.
Then $$ {\rm dim} f(A_i) = 0$$
for all $i.$ In fact, consider the general fiber $F$ of $f.$ Then $K_F = K_X \vert F$ 
and from $K_X^2 = 0$ we get $K_F^2 = 0.$ Since $\kappa (F)  = 0,$ the fiber $F$ must be
minimal, hence $K_F \equiv 0.$ Now if ${\rm dim} f(A_i) = 1$ for some $i,$ then 
$mK_X \vert F$ would be non-zero effective, contradiction. 
\sn By $K_X^2 = 0$ we obtain
$$ (\sum \lambda_i A_i)^2 = 0.$$
All $A_i$ being contained in fibers of $f,$ this implies that $\sum \lambda_i A_i \equiv \rho F$
which is impossible (i.e. $kK_X$ is generated by global sections). 

\sn (2) Now assume that $f$ is not holomorphic. Let $\pi : \hat X \La X$ be a sequence
of blow-ups such that $\hat f : \hat X \La S$ is holomorphic. Therefore we obtain some
exceptional divisor $E$ for $\hat f$ such that ${\rm dim} \hat f (E) = 1.$ It follows
that $a(\hat F) \leq 1$ for the general fiber $\hat F$ of $\hat f$, otherwise $\hat X$ 
would be algebraically connected, hence projective [Ca 81]. 
\sn Assume first that 
$a(\hat F) = 0.$ Observe $B = \pi (\hat F).$ Taking another section in $H^0(X,mK_X),$
then, restricting to $B,$ 
we obtain $D \in \vert mK_X \vert B \vert $ with $D^2 = 0.$ Hence $\pi ^*(D)^2 = 0.$
But on a (blown up) torus or K3 surface without meromorphic functions there is no such 
effective non-zero divisor. 
\sn If however $a(\hat F) = 1,$ then, arguing in the same way,
we deduce that $(\pi \vert \hat F)^*(D)$ consists of multiples of fibers of the algebraic
reduction 
$$h: \hat F \La C.$$
Since $h$ has no multi-sections, there is a map $h' : B \La C'$ to a possibly non-normal
curve $C'$ (with normalisation $C)$ such that $D$ consists of multiples of fibers of 
$h'.$ Hence $kK_X \vert B$ is generated by $H^0(B,kK_X \vert B)$ for $k \gg 0.$ In particular,
$K_X \vert B$ is nef with $\kappa = 1.$ On the other hand $K_{\hat X}Ê\vert \hat F = \pi^*(K_X \vert B) + A$
with $A$ effective. This contradicts
the fact that $\kappa (K_{\hat X}\vert \hat F) = \kappa (\hat F) = 0.$

\bn {\bf (4.17)}Ê Finally we prove Theorem 4.1   
 in the only remaining case Ê$K_{X_0}^2 \equiv 0$
but $K_{X_0}Ê\not \equiv 0$ by the same reasoning as in (4.12).

\bn Putting together (3.16) and (4.1) we now obtain
 the Main Theorem for {\it algebraically approximable} $X$ as stated in the Introduction. 

\sn Combining (3.16) with (4.10) and (4.16) we obtain as another part of the Main Theorem:

\bn {\bf 4.18 Theorem}  {\it Let $X$ be a smooth compact K\"ahler threefold with $K_X$ not nef.
Assume that $\kappa (X) = 1$ or $2$ and in case $\kappa (X) = 1 $ assume furthermore
that $X$ can be approximated algebraically. Then $X$ carries an extremal contraction.}
\vfill \eject

\noindent {\medium 5. Generic conic bundles}
\bn
{\bf (5.1)} In Part 1, (2.11) and (2.12) we had proved the following. Let $X$ be a non-projective
compact K\"ahler threefold which can be approximated algebraically or assume that $X$ is uniruled. Assume that $X$ admits a 
non-splitting 1-dimensional family $(C_t)$ of rational curves with $-K_{X_t}Ê\cdot C_t = 1,$
filling up a non-normal surface $S.$ Then the algebraic dimension $a(X) = 1,$ and $X$ has
a "generic conic bundle structure" in the following sense. 
\sn The algebraic reduction is a holomorphic map $f: X \La C$ to a smooth curve $C$ whose
general fiber is an almost homogeneous ${\bf P}_1-$bundle over an elliptic curve. The surface
$S$ is a fiber of $f.$ Moreover there exists an almost holomorphic meromorphic
map $g: X \rightharpoonup Y$ to a normal surface $Y$ such that $f$ factorises over $g$ and
such that the following holds:
{\item {(a)} if $U = f^{-1}(\{ c \in C \vert f^{-1}(c)$ is smooth or  $S = f^{-1}(c)  \},$
then $g \vert U$ is holomorphic and proper
\item {(b)} $g \vert U$ is a conic bundle 
\item {(c)}Ê$g$ is a rational quotient of $X.$} 
\sn 
The general smooth fiber of $g$ is a $\P_1$ with normal bundle $\O \oplus \O$ and therefore
we have a 2-dimensional family $(\tilde C_t)_{t \in \tilde T}$ such that for every $t_1 \in T$
there is a $t_2 \in T$ and a $t \in T$ with $C_{t_1} + C_{t_2}Ê= \tilde  C_t.$ 
\sn In general of course, the map $g$ will not be holomorphic: just perform some birational
transformation on a conic bundle. However we shall prove a criterion when $g$ is actually
holomorphic. 
\sn We say that a 2-dimensional family $(\tilde C_t)$ of rational curves splits only in the  {\bf standard
way} if the following holds.
\sn If $\tilde C_{t_0}$ is a reducible member of the 2-dimensional family, then
$\tilde C_{t_0} = C_{t_1}' + C_{t_2}'$ with smooth rational curves $C_{t_i}'$ meeting transversally
at one point and with $K_X \cdot C_{t_i}' = -1.$ 
\sn So {\it all} $\tilde C_t$ are conics. 
\sn In this terminology we have:

\bn {\bf 5.2 Theorem} {\it Let things be as in the setting (5.1) and assume that the induced
family $(\tilde C_t)$ splits only in the standard way. Then $Y$ can be taken smooth so that
$g$ is holomorphic and is a conic bundle over $Y$ whose fibers are just the $\tilde C_t.$}

\bn {\bf Proof.} Let $p: \C \La X$ be the graph of the family $(\tilde C_t)_{t \in \tilde T}$ with projection
$q: \C \La \tilde T.$ Here we consider $\tilde T$ reduced. By (5.1) there is a Zariski open
set $U \subset \tilde T$ such that $p \vert q^{-1}(U)$ is an isomorphism. Hence $p$ is
bimeromorphic, and, by considering the normalisation of $\tilde T,$ it follows that
$p$ has connected fibers. Hence $\C$ is normal and so does $\tilde T.$ It now suffices to
prove that $p$ is finite, then $p$ is biholomorphic and our claim follows. 
\sn Assume there exists $x \in X$ such that there is a curve $B \subset \tilde T$ with
$ x \in \tilde C_t$ for all $t \in B.$ Then we form the surface 
$$ A = \bigcup_{t \in B} \tilde C_t = p(q^{-1}(B)).$$ 
(a) First assume that no $\tilde C_t, t \in B$ splits. Then, $A$ being normal by Part1, (2.3),
similar arguments (even easier) as in Part1, (2.1), lead to a contradiction (by the fact that 
the $\tilde C_t$ can be moved out of $A$). 
\sn (b) So there exists $t_0 \in B$ such that $\tilde C_{t_0}$ splits. By our assumptions
the components of the splitting deform in a non-splitting family. Hence it follows from what
we have said in (5.1) (i.e. (2.12)) that there is no common point      
of the $\tilde C_t, t \in B,$ since the map $g$ is holomorphic near the reducible conics ( the
maps $g$ attached to the various 1-dimensional non-splitting families attached to the
full family of conics are of course the same). 
\sn Hence $p$ is finite, hence biholomorphic. The smoothness of $Y$ comes from deformation
theory and the smoothness of $X.$ 

\bn {\bf 5.3 Corollary} \ Ê{\it Assume in (5.1) that the $C_t$ define a geometrically extremal
ray in ${\overline {NE}}(X).$ Then $g$ is a conic bundle over a smooth surface $Y$ with
$a(Y) = 1.$}
\bn {\bf Proof.} We have to prove that the family $(\tilde C_t)$ splits only in the standard
way. Assume to the contrary that $\tilde C_{t_0} = \sum_i^p C'_i $ with $p \geq 2.$ Since 
$[C_t]$ is geometrically extremal and since $[\tilde C_{t_0}]Ê= 2[C_t]Ê\in {\overline {NE}}(X),$
we have $[C_i'] \in \overline {NE}(X),$ hence $K_X \cdot C_i' < 0.$ Hence the claim is clear.

\bn It is interesting to have a general look at conic bundles over non-algebraic surfaces.
The next proposition follows of course from (5.1) and (5.3) but it is instructive to see
a direct argument.

\bn {\bf 5.4 Proposition} {\it Let $\varphi: X \La S$ be a conic bundle over the K\"ahler
surface $S$ with $a(S) = 0$ such that $\rho (X) = \rho (S) + 1.$ Then the discriminant locus
$\Delta  = \emptyset,$ so that $\varphi$ is an analytic $\P_1-$bundle. }

\bn {\bf Proof.} Assume $\Delta \ne \emptyset.$ By the Kodaira classification $S$ is  birationally
equivalent to a torus or a K3 surface. In particular all curves in $S$ are smooth rational
curves. Now it is a basic fact on conic bundles with $\rho (X) = \rho (S) + 1$ that every
smooth rational component $C \subset \Delta $ has to meet $\overline {\Delta \setminus C}$ in
at least two points (see e.g. [Mi83], the arguments remaining true in the non-algebraic case). 
This already rules out the torus case. In the K3 case first notice that the above condition 
implies $\Delta^2 = 0.$ The union of all rational curves, and in particular $\Delta,$ is
however contractible, therefore $\Delta^2 < 0,$ contradiction. 

\bn The same type of argument together with Kodaira's classification of singular fibers of
elliptic surfaces proves

\bn {\bf 5.5 Proposition} \ Ê{\it Let $\varphi: X \La S$ be a conic bundle over a K\"ahler surface
$S$ with $a(S) = 1,$ satisfying $\rho (X) = \rho (S) + 1.$ Let $\Delta$ be its discrminant
locus. Then $\Delta = \bigcup F_i,$ where the $F_i$ are smooth fibers of the algebraic 
reduction $f: S \La C$ (which is an elliptic fibration) or reduction of multiple smooth
fibers.}

\bn {\bf 5.6 Proposition} {\it Let us assume the situation of (5.5). Then $a(X) = 1$ and 
$f \circ \varphi$ is an algebraic reduction of $X.$}

\bn {\bf Proof.} If $a(X) \ne 1,$ then $a(X) = 2.$ This case is already ruled out in Part 1.
The argument is as follows. Let $g: X \rightharpoonup Z$ be an
algebraic reduction and $F$ a general fiber of $\varphi.$ Since $g$ is an "elliptic fibration",
we have ${\rm dim} g(F) = 1.$ It follows easily that any two points in $X$ can be
 joined by a chain of curves, i.e. $X$ is algebraically connected. Hence $X$ is Moishezon by
[Ca81], contradiction.

\vskip.5cm \noindent
{\medium 6. Minimal Models}
\bn {\bf (6.1)} The weak minimal model conjecture (WMMC) in the K\"ahler case predicts
that every compact K\"ahler manifold $X$ is either uniruled or birationally equivalent
to a K\"ahler $n-$fold $X'$ with at most terminal singularities such that $K_{X'}$ is
nef; such an $X'$ is called minimal model.
\sn The abundance conjecture (AC) says that every minimal model is semi-ample, i.e.
some multiple $mK_{X'}$ is generated by global sections. A minimal model with 
semi-ample $K_{X'}$ is also called good minimal model. 
\sn The strong minimal model conjecture (SMMC) states that starting from a compact K\"ahler
manifold $X$ one get derive either a ${\bf Q}-$Fano fibration or a minimal model by
a sequence of birational divisorial contractions or flips.
 
\bn See [KMM87] for the background in the algebraic case. SMMC holds for algebraic threefolds
by [Mo88] and AC by [Mi88] and [Ka92]. 

\bn {\bf (6.2)} Both WMMC and AC have been proved by Nakayama [Na88] in case the compact
K\"ahler threefold carries an elliptic fibration, in particular if $a(X) = 2.$   
   
\bn {\bf (6.3)} We next describe the structure theorem of Fujiki [Fu83] for non-algebraic
non-uniruled
compact K\"ahler threefolds with $a(X) \leq 1.$
\sn (a) If $a(X) = 1,$ then either we have a holomorphic algebraic reduction $f: X \La C$
with the general smooth fiber being a torus or $X$ is birationally $(C \times S)/G$     
with $C$ a compact Riemann surface, $S$ a torus or a K3-surface with $a(S) = 0$ and $G$ a 
finite group acting on both $C$ and $S$ and on $C \times S$ by $g(x,y) = (gx,gy).$ 
\sn (b) If $a(X) = 0, $ then either $X$ is a Kummer manifold, i.e. $X$ is birationally
$T/G,$ where $T$ is a torus and $G$ a finite group or $X$ is simple, i.e. 
$X$ does not carry a covering family of compact subvarieties and does not admit a
meromorphic map to a Kummer manifold.

\bn {\bf 6.4 Theorem} {\it Kummer threefolds with $a(X) = 0$ have good minimal models.}

\bn {\bf Proof.} Let $X$ be a Kummer threefold with $a(X) = 0.$ So $X$ is birationally
$Y = T/G$ with a torus $T$ of algebraic dimension $0.$ Since $T$ does not carry positive-dimensional 
subvarieties, $Y$ can have only isolated singularities and these are quotient singularities,
since $g: T \La T/G$ is unramified outside a finite set. In particular $Y$ has only 
canonical singularities (see e.g. [KMM87]). By Reid [Re83] there is a partial resolution
$$ f: X' \La Y $$
such that $X'$ has terminal singularities and $f$ is crepant, i.e. $K_{X'}Ê= f^*(K_Y).$
Since obviously $K_Y \equiv 0, X'$ is a minimal model. Since $K_T = g^*(K_Y),$ we 
even have $mK_Y = \O_Y,$ hence $X'$ is good. 

\bn {\bf (6.5) Remark} (1) It might also be possible to construct directly a minimal
model in case $X \sim (C \times S)/G,$ but the case of the holomorphic torus fibration
is certainly harder. Possibly one has to prove SMMC in that case to proceed to a minimal
model. See (6.7) for more comments on SMMC. 
\sn (2) A consequence of the existence of good minimal models is the non-existence of
simple threefolds as in (6.3) defined; see the Motiviation following the Introduction.

\bn We finally show that once we know the existence of good minimal models, then the notions
of algebraic nefness and nefness coincide (for threefolds)

\bn {\bf 6.6 Theorem} {\it Let $X$ be a smooth compact K\"ahler threefold. Assume that $X$
has a good minimal model. Then $K_X$ is nef if and only if $K_X \cdot C \geq 0$ for
all curves $C \subset X.$} 

\bn {\bf Proof.} Due to (4.10) and (4.15) we could restrict ourselves to the case $\kappa (X) = 0.$
However we will give a simultaneous proof in all cases making the arguments independent of
section 4.
\sn One direction being obvious we assume that $K_X \cdot C \geq 0$ for all curves $C \subset X.$
Let $$ h: X \rightharpoonup X'$$ be a bimeromorphic map to a good minimal model. 
Choose a sequence of blow-ups $f: \hat X \La X$ such that the induced
map $g: \hat X \La X'$ is a morphism. Let
$$ D' \in \vert mK_{X'}Ê\vert $$ be a general smooth member and $D \subset X$ be its strict
transform (in case $\kappa (X) = 0,$ we have $ D' = 0)$. Then we can write 
$$ mK_X = D + \sum \mu_i A_i, \eqno (*)$$
with $mu_i > 0.$ To determine the structure of the $A_i,$ we write
$$ K_{\hat X} = g^*(K_{X'}) + \sum a_j F_j, a_j > 0.$$ 
Let $\hat D$ be the strict transform of $D$ in $\hat X.$ Then 
$$ g^*(D') + \sum ma_j F_j = \hat D + \sum c_j F_j \in \vert mK_{\hat X}Ê\vert$$
with $c_j \geq 0.$ 
Therefore 
$$D + \sum c_j f(F_j)_0 \in \vert mK_X \vert,$$
where the index $0$ indicates that $f(F_j)$ is omitted if ${\rm dim} f(E_j) \leq 1.$
So $A_i = f(E_j)_0$ and in particular all $A_i$ are algebraic. 
\sn By (4.9) it is sufficient to show that $K_X \vert D$ is nef and that $K_X \vert A_i$
is nef for all $i.$ Since the $A_i$ are algebriac, the second statement is clear. 
\sn (1) If $\kappa (X) = 0,$ then $D' = 0, $ hence $D = 0$ and there is nothing to prove.
\sn (2) If $\kappa (X) \geq 1,$ we note that $mK_X \vert D$ is effective, hence we conclude by 
(4.14).

\bn {\bf 6.7 Corollary}  {\it Assume that WMMC and AC hold in dimension 3. Let $X$ be a compact
K\"ahler threefold with at most terminal singularities. Then $K_X$ is nef if and only
$K_X \cdot C \geq 0$ for all curve $C \subset X.$}

\bn {\bf 6.8 Remark} The most difficult step in the construction of minimal models for
algebraic threefolds is of course the existence of flips. However this is reduced in [Ka88]
to a local analytic problem around the rational curves which have to be flipped. [Ka88]
shows also that in the analytic this local reduction works. Since Mori proves in [Mo88]
the local analytic existence of flips, we have already the existence of flips in the 
analytic category. The proof of termination is just the same as in the algebraic case.

\bn \bn \bn {\bf References}
\bn \bn
{\item {[AB93]} Alessandrini,L;Bassanelli,G.: Plurisubharmonic currents and their extensions
across analytic subsets. Forum Math. 5, 577 - 602 (1993)

\item {[Al96]} Alessandrini,L.: Letter to the author, february 1996

\item {[AP]} Andreatta.M.; Peternell,T.: On the limits of manifolds with nef canonical bundles.
To appear in Trento Proc. Conference Algebraic Geometry 1995

\item {[Ba94]} Bassanelli,G.: A cut-off theorem for plurisubharmonic currents. Forum
Math. 6, 567 - 595 (1994)

\item {[Be83]} Beauville,A.: Vari\'et\'es k\"ahleriennes dont la premiere classe de Chern
est nulle. J. Diff. Geom. 18, 755 - 782 (1983)

\item {[BPV84]} Barth,W.;Peters,C.;van de Ven,A.: Compact complex surfaces. Erg. d. Math., 3.
Folge, Band 4. Springer 1984

\item {[Ca81]} Campana,F.: Cor\'eduction alg\'ebrique d'un espace analytique faiblement
k\"ahl\'erien compact. Inv, math. 63, 187 - 223 (1981)

\item {[CP94]} Campana,F;Peternell,T.: Towards a Mori theory on compact K\"ahler threefolds,I.
To appear in Math. Nachr. 1996

\item {[De92]} Demailly,J.P.: Regularisation of closed positive currents and intersection
theory. J. Alg. Geom. 1, 361 - 410 (1992) 

\item {[DPS94]} Demailly,J.P.;Peternell,T.;Schneider,M.: Compact complex manifolds with
numerically effective tangent bundles. J. Alg. Geom. 3, 295 - 345 (1994) 

\item {[Fu83]} Fujiki,A.: On the strucutre of compact complex manifolds in class ${\cal C}.$
Adv. Stud. Pure Math. 1, 231 - 302 (1983)

\item {[Fj78]}ÊFujita,T.: K\"ahler fiber spaces over curves. J. Math. Soc. Japan 30, 779 - 794
(1978)

\item {[HL83]} Harvey,R.;Lawson,B.: An intrinsic characterisation of K\"ahler manifolds. 
Inv. math. 74, 169 - 198 (1983)

\item {[Ka81]} Kawamata,Y.: Characterisation of abelian varieties. Comp. math. 43, 253 - 276
(1981)

\item {[Ka85]} Kawamata,Y.: Minimal models and the Kodaira dimension of algebraic
fiber spaces. J. reine u. angew. Math. 363, 1 - 46 (1985)

\item {[Ka86]} Kawamata,Y.: On the plurigenera of minimal algebraic threefolds with
$K \equiv 0.$ Math. Ann. 275, 539 - 546 (1986)

\item {[Ka88]} Kawamata,Y.: The crepant blowing-ups of 3-dimensional canonical singularities
and its application to degeneration of surface. Ann. Math. 127, 93 - 163 (1988) 

\item {[Ka92]} Kawamata,Y.: Abundance theorem for minimal threefolds. Inv. math. 108,
229 - 246 (1992)

\item {[KMM87]} Kawamata,Y.;Matsuda,K.;Matsuki,K.: Introduction to the minimal model problem.
Adv. Stud. Pure Math. 10, 283 - 360 (1987)

\item {[Ko91]} Koll\'ar, J.: Extremal rays on smooth threefolds. Ann. scient. Ec. Norm. Sup.
24, 339 - 361 (1991)

\item {[Mi87]} Miyanishi,M.: Algebraic threefolds. Adv. Stud. Pure Math. 1, 69 - 99 (1983)

\item {[Mi87]} Miyaoka,Y.: The Chern classes and Kodaira dimension of a minimal
variety. Adv. Stud. Pire Math. 10, 449 - 476 (1987)  
Adv. Stud. Pure Math. 10, 

\item {[Mi88]} Miyaoka,Y.: Abundance conjecture for threefolds: $\nu = 1$ case. Comp. math. 68,
203 - 220 (1988)

\item {[Mo82]} Mori,S.: Threefolds whose canonical bundles are not numerically effective.
Ann. Math. 116, 133 - 176 (1982)

\item {[Mo88]} Mori,S.: Flip theorem and the existence of minimal models for 3-folds.
J. Amer. Math. Soc. 1, 117 - 253 (1988)

\item {[Na88]} Nakayama,N.: On Weierstrass models. Alg.Geom. and Comm. Algebra; vol. in
honour of Nagata, vol. 2, 405 - 431. Kinokuniya, Tokyo 1988

\item {[Pe95]} Peternell,T.: Minimal varieties with trivial canonical class,I. Math. Z.
217, 377 - 407 (1994)                             

\item {[Re83]} Reid,M: Minimal models of canonical threefolds. Adv. Stud. Pure Math. 1,
131 - 180 (1983)

\item {[Sk82]} Skoda,H.: Prolongement des courants, positifs, ferm\'es, de masse finie.
Inv. math. 66, 361 - 376 (1982)

\item {[Su88]} Sugiyama,Y.: A geometry of K\"ahler cones. Math. Ann. 281, 135 - 144 (1988)

\item {[Ue87]} Ueno,K.: On compact analytic threefolds with non-trivial Albanese torus.
Math. Ann. 278, 41 - 70 (1987)

\item {[Wi89]} Wisniewski,J.: Length of extremal rays and generalised adjunction. Math. Z.
200, 409 - 427 (1989) }

\bn \bn \bn 

\noindent Thomas Peternell, Math. Institut der Univ. Bayreuth, 
\sn D-95440 Bayreuth, Germany

\end